\documentclass{article}

\usepackage{fullpage}

\usepackage{amsmath}
\usepackage{amssymb}
\usepackage{amsthm}
\usepackage{mathabx}
\usepackage{mathrsfs}
\usepackage{tabularx}
\allowdisplaybreaks

\usepackage{cancel}

\usepackage{color}
\usepackage[colorlinks=true, pdfstartview=FitV, linkcolor=blue, citecolor=blue, urlcolor=blue,pagebackref=false]{hyperref}

\let\emptyset \undefined
\newsymbol\emptyset    203F

\newcommand{\ud}{\operatorname{d}\! }
\newcommand{\R}{\mathbb{R}}

\newcommand{\N}{\mathbb{N}}
\newcommand{\hol}{H\"older }

\newcommand{\V}{\mathbb{V}}

\def \E{\mathbb E}
\def \Eb{\mathbf E}
\def \P{\mathbb P}
\def \Pb{\mathbf P}

\def \H{\mathscr H}

\def \C{\mathscr C}
\def \L{\mathcal L}

\def \Tps{\mathscr T} 
\def \S{\mathcal S}

\def \Nop{\mathcal N}

\def \A+{\mathscr A_+}
\def \cpartial {\partial^\perp}
\def \flap {(-\Delta)}
\def \cGFF {\underline{\omega}}
\def \V {\mathcal V}
\def \Vs {\mathbf V}
\def \S {\mathcal S}
\newcommand{\lb}{\langle}
\newcommand{\rb}{\rangle}
\newcommand{\vertiii}[1]{{\left\vert\kern-0.25ex\left\vert\kern-0.25ex\left\vert #1 
		\right\vert\kern-0.25ex\right\vert\kern-0.25ex\right\vert}}

\newtheorem{theorem}{Theorem}

\newtheorem{lemma}[theorem]{Lemma}
\newtheorem{proposition}[theorem]{Proposition}
\newtheorem{remark}[theorem]{Remark}

\begin{document}
	
	\title{Brownian particle in the curl of 2-d stochastic heat equations}
	\author{Guilherme L. Feltes$^1$ and Hendrik Weber$^2$}
	\date{}
	\maketitle

\begin{center}
	{\footnotesize Department of Mathematical Sciences,  
		
		University of Bath
	
and 

Institute for Analysis and Numerics,

Universit\"at M\"unster.}
\end{center}
\footnotetext[1]{\texttt{guilherme.feltes@uni-muenster.de}}
\footnotetext[2]{\texttt{hendrik.weber@uni-muenster.de}}	

\begin{abstract}
	\noindent  We study the  long time behaviour of a Brownian particle evolving in a  dynamic random environment. Recently, Toninelli et al.  \cite{Cannizaroetal2021} proved sharp $\sqrt{log}$-super diffusive bounds for a Brownian particle in the curl of (a regularisation of) the 2-d Gaussian Free Field (GFF) $\cGFF$. We consider a  one parameter family of Markovian and Gaussian dynamic environments which are reversible with respect to the law of $\cGFF$. Adapting their method, we show that if $s\ge1$, with $s=1$ corresponding to the standard stochastic heat equation, 
	then the particle stays $\sqrt{log}$-super diffusive, whereas if $s<1$,
	corresponding to a fractional heat equation,
	then the particle becomes diffusive. In fact, for $s<1$, we show that this is a particular case of Komorowski and Olla \cite{komorowskiandolla2003}, which yields an invariance principle through a Sector Condition result. Our main results agree with the Alder-Wainwright scaling argument (see \cite{WainwrightandAlder1967,WainwrightandAlder1970,WAG1971,Fosteretal1977}) used originally in \cite{TothandValko2012} to predict the $\log$-corrections to diffusivity.
	We also provide examples 
	which display $log^a$-super diffusive behaviour for $a\in(0,1/2]$.
	
	
\end{abstract}

\setcounter{tocdepth}{3}

\section{Introduction and main result} \label{section_introduction}

We study the motion of a Brownian particle in $\R^2$, evolving in a dynamic random environment (DRE), given by the solution to the It\^o SDE
\begin{equation} \label{intro_def_main_SDE}
	\begin{cases}
		& \ud X(t) =  \omega_t(X(t)) \ud t  + \sqrt{2}\ud B(t), \quad t \geq 0, \\
		& X(0) = 0 \; ,
	\end{cases}
\end{equation}
where $(B(t))_{t \ge 0}$ is a standard two-dimensional Brownian motion and $(\omega_t(x))_{t\ge0,x\in\R^2}$ is a time-dependent  random field which is independent from $(B(t))_{t \ge 0}$. 
We take $(\omega_t(x))_{t\ge0,x\in\R^2}$ to be a regularised version of the curl of the solution to the (fractional) stochastic heat equation with additive noise in $\R^2$ and initial condition given by the curl of the regularised Gaussian Free Field (GFF) $\cGFF$. The coordinates of $\omega_t = (\omega_t^1,\omega_t^2)$ satisfy
\begin{equation} \label{intro_def_environmental_process_omega}
	\begin{cases}
		\ud \omega_t^k = -\flap^s\omega_t^k\ud t + \sqrt{2} \cpartial_k \flap^{\frac{s-1}{2}}\ud W_t \; , \; t \ge 0 \; \; , \; \; k =1,2 \; \; ,\\
		\omega_0 = \cGFF \; \; ,
	\end{cases}
\end{equation}
where $s \in [0,\infty)$ and $(\cpartial_1,\cpartial_2) := (\partial_{x_2},-\partial_{x_1})$. Here, $W$ is a mollified (in space) space-time white noise, with covariance $\Eb[W_r(x)W_t(y)] = \min\{r,t\} V(x-y)$, and $\cGFF$ is distributed according to the law of the curl of a mollified GFF. More precisely, for every $k,l = 1,2$, $r,t \ge 0$ and $x,y \in \R^2$, $W_t^s(x) := \flap^{\frac{s-1}{2}}W_t(x)$ and $\cGFF$ have mean zero and covariance 
\begin{align}
	& \Eb[\cpartial_kW_r^s(x)\cpartial_lW_t^s(y)] = \min\{r,t\}\cpartial_k\cpartial_lV \convolution g_{1-s}(x-y) = \min\{r,t\}(\cpartial_k(-\Delta)^{s-1}\delta_x,\cpartial_l\delta_y)_V \text{ if } s \le 1 \; , \label{intro_def_covariance_curl_mollified_transf_white_noise}\\
	& \Eb[\cpartial_kW_r^s(x)\cpartial_lW_t^s(y)] = \min\{r,t\}\cpartial_k\cpartial_l(-\Delta)^{s-1}V(x-y) = \min\{r,t\}(\cpartial_k(-\Delta)^{s-1}\delta_x,\cpartial_l\delta_y)_V  \text{ if } s > 1 \; , \label{intro_def_covariance_curl_mollified_transf_white_noise_s>1}\\
	& \Eb[\cGFF^k(x)\cGFF^l(y)] = -\cpartial_k\cpartial_lV \convolution g_1(x-y) = (\cpartial_k(-\Delta)^{-1}\delta_x,\cpartial_l\delta_y)_V \; , \label{intro_def_covariance_curl_mollified_GFF}
\end{align}
where $\convolution$ denotes the convolution over $\R^2$ and, for every $\varphi_1,\varphi_2 \in \S(\R^2)$, the space of Schwartz functions over $\R^2$,
\begin{equation} \label{intro_def_mollified_inner_product}
	(\varphi_1,\varphi_2)_V := \int_{\R^2}\int_{\R^2}\varphi_1(x)V(x-y)\varphi_2(y)\ud x\ud y \; \; .
\end{equation}
The smooth function $V$ is given by
\begin{equation} \label{intro_def_mollifier_V}
	V := U \convolution U \; ,
\end{equation}
for a $U \in C^\infty(\R^2)$, radially symmetric, decaying exponentially fast at infinity and with $\int_{\R^2} U(x) \ud x = 1$. To simplify some computations, we may also assume that $U$ has Fourier transform supported in ball of radius $1$. Also, the kernel $g_r : \R^2 \setminus \{0\} \to \R$ in \eqref{intro_def_covariance_curl_mollified_transf_white_noise} and \eqref{intro_def_covariance_curl_mollified_GFF} is given by
\begin{equation*}
	\begin{cases}
		g_{1}(x) &= - (2\pi)^{-1}\log |x|  \; \; ; \\
		g_{r}(x) &= \dfrac{\Gamma(1-r)}{4^r\Gamma(r)\pi} \dfrac{1}{|x|^{2-2r}} \quad \text{ if } r \in (0,1) \; \; ; \\
		g_{0} & = \delta_0 \;  \; \; ,
	\end{cases}
\end{equation*}
where $\Gamma$ denotes the Gamma function. I.e., $g_r$ is the Green's function of $(-\Delta)^r$ in $\R^2$, for $r \in [0,1]$. Also, the fractional Laplacian $(-\Delta)^{s-1}$ for $s > 1$ can be defined in terms of its Fourier multiplier, as $\widehat{(-\Delta)^{s-1}f}(p) = |p|^{2(s-1)} \widehat{f}(p)$.

\begin{remark} \label{intro_remark_mollifier_makes_everything_well_defined}
	Note that expressions \eqref{intro_def_covariance_curl_mollified_transf_white_noise}-\eqref{intro_def_covariance_curl_mollified_GFF} make sense due to the presence of the smooth function $V$. Plugging, e.g., the right-hand side of \eqref{intro_def_covariance_curl_mollified_GFF} into \eqref{intro_def_mollified_inner_product}, we get
	\begin{equation*}
		(\cpartial_k(-\Delta)^{-1}\delta_x,\cpartial_l\delta_y)_V =
		\int_{\R^2}\int_{\R^2}\delta_x(u) \cpartial_k\cpartial_l(-\Delta)^{-1}V(u-v)\delta_y(v)\ud u\ud v \; ,
	\end{equation*}
	which is equal to the expression in the middle of \eqref{intro_def_covariance_curl_mollified_GFF}. Furthermore, even though 
	the GFF in the full space is only defined up to a constant (i.e. inverting the Laplacian $\Delta$), taking the derivatives $\cpartial_k\cpartial_l$ of its regularisation makes it rigorous without ambiguity. The same reasoning holds to define the noise $\cpartial_k\flap^{\frac{s-1}{2}}W_t$ when $s=0$.
	
\end{remark}

\begin{remark}
	As we show in  Proposition \ref{op_fock_prop_invariant_measure_fSHE}, the dynamics \eqref{intro_def_environmental_process_omega} leave the law of  $\cGFF$ invariant.
	 The case $s=1$ corresponds to the standard stochastic heat equation (SHE), whereas $s=0$ is the infinite dimensional Ornstein Uhlenbeck process, as defined, for example in \cite[Chapter 1.4]{Nualart2006}.
The parameter $s \in [0,\infty)$ controls the speed of the environment on different scales: smaller values of $s$ correspond to faster movement of the larger scales. 

\end{remark}


By definition, the drift field $\omega_t(x)$ in \eqref{intro_def_main_SDE} is \textit{divergence-free}. Brownian particles evolving in stationary divergence-free random fields  have been considered as a toy model for anomalous diffusions in inhomogeneous media, such as the motion of a tracer particle in an incompressible turbulent flow. See e.g. the surveys \cite[Chapter 11]{landimetal2012} and \cite{Kozlov1985}.
Depending on the decay of the spatial correlations of the drift field, the particle could behave either \textit{diffusively} or \textit{superdiffusively}, meaning that the mean square displacement satisfies for large $t$
\begin{equation} \label{intro_def_diffusive_vs_superdiffusive}
D(t) := \frac{\E[|X(t)|^2]}{t} \approx
	\begin{cases}
	 1 \quad \quad \quad \quad \quad \text{diffusive,} \\
	 t^\nu  \; , \; \nu > 0 \; \quad \text{superdiffusive.}
	\end{cases}
\end{equation}
Here, $\E$ denotes the expectation under the joint law of $B$ and  $\omega$, see Section \ref{section_set_main_result}. If the correlations of the environment decay fast enough (see e.g. \cite[Chapter 11]{landimetal2012}), one gets  diffusive behaviour, and if the decay is too slow (see \cite{KomorowskiandOlla2002}), one gets  superdiffusive behaviour.  There is, however, an intermediate regime for which the correlations decay in such a way that $D(t)$ diverges only as $(\log t)^\gamma$, for $\gamma > 0$. These logarithmic corrections 
are expected to be present in two-dimensional Brownian particles evolving in isotropic random drift fields. Indeed, by the Alder-Wainwright scaling argument (see \cite{WainwrightandAlder1967,WainwrightandAlder1970,WAG1971,Fosteretal1977}), in 2d, if the displacement of the particle scales faster than the correlations of the environment field, then the (only) expected behaviour for the mean square displacement of the particle is to be of order $t \sqrt{\log t}$. We briefly elaborate on this, following the Appendix of \cite{TothandValko2012}. Let $K(t,x) := \Eb[\omega_0(0)\omega_t(x)]$. Now, assume that $\P(X(t) \in \ud x) \approx \alpha(t)^{-2} \varphi(\alpha(t)^{-1}x) \ud x$, where $\varphi$ is a density and $\alpha(t) = t^\nu (\log t)^\gamma$ for some $\nu, \gamma \ge 0$. If we also assume that $K(t,x) \approx \beta(t)^{-2}\psi(\beta(t)^{-1}x)$, for another density $\psi$,
then if 
\begin{equation} \label{intro_scaling_argument_hip}
	\dfrac{\beta(t)}{\alpha(t)} \le \text{constant, for } t \ge 0 \; \; ,
\end{equation}
we must have $\nu = 1/2$ and $\gamma = 1/4$, which yields $X(t) \approx t \sqrt{\log t}$. We emphasise here that this argument, even though instructive, it is not mathematically rigorous.
Indeed, the $\sqrt{\log}$ correction was rigorously established recently by Toninelli et al.  \cite{Cannizaroetal2021}. They showed that for a time-independent drift field $\omega$  distributed according to the law of $\cGFF$, one has
\begin{equation} \label{intro_sharp_estimates_CHT21}
	D(t) \approx \sqrt{\log t} \text{ as } t \to \infty \; \; ,
\end{equation}
up to $\log \log t$ corrections, confirming a conjecture made by T{\'o}th and Valk{\'o} \cite{TothandValko2012} based on this scaling argument.
The result
was obtained in the \textit{Tauberian sense}\footnote{For a discussion on the connection between the asymptotics of $D(t)$ and $D_T(\lambda)$, see \cite[Lemma 1]{Tothetal2012} or \cite[Remark 2.3]{Cannizaroetal2021}}, i.e., in terms of the Laplace transform of the mean square displacement
\begin{equation} \label{set_main_res_def_D}
	D_T(\lambda) := \int_0^\infty e^{-\lambda t} \E[|X(t)|^2] \ud t \; , \; \; \lambda > 0 \; .
\end{equation}

Note that in the case considered by \cite{TothandValko2012,Cannizaroetal2021}, the correlations of the drift field do not scale in time since the drift field is time-independent, so \eqref{intro_scaling_argument_hip} is trivially satisfied.
Moving to the time-dependent case treated in the present work, if we take $s\ge1$ in \eqref{intro_def_environmental_process_omega}, we still have that the correlations of the drift field do not scale fast enough - $t^{\frac{1}{2s}}$ for $\omega$ vs. $t^{\frac{1}{2}} (\log t)^{\frac{1}{4}}$ for $X$. Therefore, we should still expect for the particle $X$ to behave \textit{$\sqrt{\log}$-superdiffusively}, since condition \eqref{intro_scaling_argument_hip} remains true. However, if we move to the case where $s < 1$ in \eqref{intro_def_environmental_process_omega}, then the picture changes substantially and condition \eqref{intro_scaling_argument_hip} is no longer satisfied, since $t^{\frac{1}{2}} (\log t)^{\frac{1}{4}} << t^{\frac{1}{2s}}$ for $t >> 1$. Theorem \ref{intro_theo_main_theo} and Theorem \ref{intro_theo_invariance_principle} below rigorously establish the expected abrupt difference between super diffusive and diffusive behaviours depending on the exponent $s$, agreeing with the scaling argument.

\begin{theorem} \label{intro_theo_main_theo}
	If $s \ge 1$ in \eqref{intro_def_environmental_process_omega}, then, for every $\varepsilon > 0$, there exist constants $A_\varepsilon, B_\varepsilon > 0$, depending only on $\varepsilon$ and $s$, such that, for $\lambda \in (0,1)$, we have
		\begin{equation} \label{intro_super_dif_bounds_s>1}
			A_\varepsilon(\log |\log \lambda|)^{-1 - \varepsilon}
			\le \lambda^2 \dfrac{D_T(\lambda)}{\sqrt{|\log \lambda|}} \le
			B_\varepsilon(\log |\log \lambda|)^{1 + \varepsilon} \; .
		\end{equation}
\end{theorem}

\vspace{3mm}

\noindent For the case $s \in [0,1)$ we can apply a sector condition result of Komorowski and Olla \cite{komorowskiandolla2003} to obtain the following invariance principle.

\begin{theorem} \label{intro_theo_invariance_principle}
	If $s \in [0,1)$ in \eqref{intro_def_environmental_process_omega}, then there exist constants $A, B > 0$, such that, for all $t \ge 0$, we have
	\begin{equation} \label{intro_dif_bounds_s<1}
		A
		\le D(t) 
		\le
		B \; .
	\end{equation}
	Furthermore, let $(Q^\omega_\varepsilon)_{\varepsilon\in (0,1]}$ denote the laws of $(\varepsilon X(\frac{t}{\varepsilon^2}))_{t\ge0}$, over $C[0,\infty)$, for $\varepsilon \in (0,1]$, given the initial configuration $\omega_0 = \omega$. Then $(Q^\omega_\varepsilon)_{\varepsilon\in (0,1]}$ converge weakly, with respect to the law of $\cGFF$,
	as $\varepsilon \downarrow 0$, to the law of a Brownian motion with deterministic covariance matrix $D$, which only depends on $s$. The covariance matrix $D$ is defined in \eqref{proof_theo_inv_principle_covariance_matrix_D}.
\end{theorem}

\vspace{3mm}

The asymptotic behaviour of $D_T(\lambda)$ in \eqref{intro_super_dif_bounds_s>1} is a reflection of the fact that the dynamics provided by the SHE (with the full Laplacian) does not mix the environment fast enough to produce a scaling of the correlations which is faster than the scaling of the displacement of the particle, as discussed above. On the other hand, the result in \eqref{intro_dif_bounds_s<1} confirms that the fractional dynamics on the environment changes dramatically the behaviour of the particle.
Moreover, the estimates in \eqref{intro_super_dif_bounds_s>1} are exactly the same as the ones obtained in \cite{Cannizaroetal2021}, and our proof is an adaptation of theirs, which is based on Yau's method \cite{Yau2004} of recursive estimates of iterative truncations of the resolvent equation in \eqref{set_main_res_laplace_transform_to_resolvent_eq}. Indeed, when $s\ge1$, the dominant terms in the estimates are the ones coming from the stationary drift field, which are the same as for the static case. What we show is that we can remove the additional terms coming from the dynamics of the environment in the estimates, maintaining the same asymptotic behaviour. However, when $s<1$, the dominant terms are now precisely the ones coming from the dynamics of the environment. The effect can be seen already in the first upper bound obtained by the first truncation of \eqref{set_main_res_laplace_transform_to_resolvent_eq}, and it is enough to show \eqref{intro_dif_bounds_s<1} in Theorem \ref{intro_theo_invariance_principle}, see Remark \ref{anal_res_remark_level_1}. 



If now we consider intermediate regimes between $s=1$ and $s<1$, only adding a logarithmic divergence to the operator $\Delta$ in \eqref{intro_def_environmental_process_omega},
we obtain something which was not predicted by the Alder-Wainwright scaling argument. Namely, for any given $a \in (0,\frac{1}{2}]$, we can find an interpolation between the regimes $s=1$ and $s<1$ such that we prove corrections to diffusivity of order $(\log t)^a$.
More precisely, if we consider that the coordinates of $\omega_t = (\omega_t^1,\omega_t^2)$ satisfy
\begin{equation} \label{intro_def_environmental_process_omega_log_div}
	\begin{cases}
		\ud \omega_t^k = (\log(e + (-\Delta)^{-1}))^\gamma\Delta\omega_t^k\ud t + \sqrt{2} (\log(e + (-\Delta)^{-1}))^\frac{\gamma}{2}\cpartial_k \ud W_t \; , \; t \ge 0 \; \; , \; \; k =1,2 \; \; ,\\
		\omega_0 = \cGFF \; \; ,
	\end{cases}
\end{equation}
for a parameter $\gamma > 0$. Then, we can show the following
\begin{theorem} \label{intro_theo_log_div}
	If $(X(t))_{t\ge0}$ is the solution to \eqref{intro_def_main_SDE} with $(\omega_t)_{t\ge0}$ solution to \eqref{intro_def_environmental_process_omega_log_div}, then, for every $\gamma \in [\frac{1}{2},\infty)$, there exist constants $A,B>0$, only depending on $\gamma$, such that:
	
	\noindent If $\gamma \in [\frac{1}{2},1)$, then for $\lambda \in (0,1)$,
	\begin{equation} \label{intro_loglog_bounds}
		A \le \dfrac{\lambda^2 D_T(\lambda)}{|\log \lambda|^{1-\gamma}} \le B \; .
	\end{equation}
	\noindent If $\gamma = 1$, then for $\lambda \in (0,1)$,
	\begin{equation} \label{intro_loglog_bounds_gamma_1}
		A \le \dfrac{\lambda^2 D_T(\lambda)}{\log |\log \lambda|} \le B \; .
	\end{equation}
	Furthermore, if $\gamma > 1$, we have
	\begin{equation} \label{intro_dif_bounds_gamma>1}
		A
		\le D(t) 
		\le
		B \; .
	\end{equation}
\end{theorem}

\subsection{Structure of the paper}
In Section \ref{section_set_main_result} we define the 
environment seen from the particle process as a technical tool. In Section \ref{section_operators_fock_space} we derive the action of the infinitesimal generator of the environment seen from the particle on Fock space, and show that the law of $\cGFF$ is invariant under the family of dynamics given by \eqref{intro_def_environmental_process_omega}. Section \ref{section_iterative_anal_resolvent_eq} contains the proof of the main recursive estimates through an iterative analysis of the resolvent equation in \eqref{set_main_res_laplace_transform_to_resolvent_eq} and a proof of \eqref{intro_dif_bounds_s<1} in Theorem \ref{intro_theo_invariance_principle} using only the first truncation of the resolvent equation. In Section \ref{section_proof_main_theo_diffusive_bounds} we prove Theorem \ref{intro_theo_main_theo} by using the recursive estimates obtained in Section \ref{section_iterative_anal_resolvent_eq}. In Section \ref{section_proof_theo_inv_princ}, we present a general overview of the method in \cite{komorowskiandolla2003} of homogenisation of diffusions in divergence-free, Gaussian and Markovian fields and show that for $s<1$ we may apply their results to get Theorem \ref{intro_theo_invariance_principle}. In Section \ref{section_proof_theo_log_div} we prove Theorem \ref{intro_theo_log_div}. Appendices \ref{section_appendix_their_lemmas} and \ref{section_appendix_our_lemmas} gather important ingredients from Toninelli et al. \cite{Cannizaroetal2021}, and some generalisations to the present setting,
necessary in Sections \ref{section_iterative_anal_resolvent_eq} and \ref{section_proof_main_theo_diffusive_bounds}
and Appendix \ref{section_appendix_end_proof_theo_1} presents the final argument for the proof of Theorem \ref{intro_theo_main_theo}, taken from \cite{Cannizaroetal2021}.

\section{Setting and preliminaries} \label{section_set_main_result}

Let $\Tps_0 := (\Omega, \mathcal B, \Pb)$ be a probability space supporting $\cGFF$ and an independent Wiener process $W$ as defined between displays \eqref{intro_def_environmental_process_omega} and \eqref{intro_def_covariance_curl_mollified_transf_white_noise}. 
Let $\Tps_1 := (\Sigma, \mathcal F, Q)$ be another probability space supporting a standard 2d Brownian motion $B$.  We consider solutions to the system \eqref{intro_def_main_SDE}, \eqref{intro_def_environmental_process_omega} on $\Omega \times \Sigma$ equipped with the product measure  $\P = \Pb \otimes Q$. The law of $(X(t))_{t\ge0}$ under $\P$ is called the  \textit{annealed} law. Note, that under $\P$, the process $(X(t))_{t\ge0}$ alone  is not Markovian. Notwithstanding, we may define a different Markovian process,  the so-called \textit{environment seen from the particle}, which takes values on the larger space of functions over $\R^2$ \cite{Komorowski2000}. It evolves by spatially shifting the environment by the position of the walker, at any given time $t\ge 0$. Precisely, we set
\begin{equation} \label{set_main_res_def_eta_t}
	\eta_t :
	 = \omega_t(\cdot + X(t)) \; , \; t\ge 0 \; \; .
\end{equation}
The law of $X$  is rotationally invariant, and therefore we have that $\E[|X(t)|^2] = \E[X_1(t)^2 + X_2(t)^2] = 2\E[X_1(t)^2]$. Hence we
 may focus on its first coordinate only. Furthermore, $\E[X(t)] = 0$. Formula \eqref{set_main_res_def_eta_t}  allows us to write
\begin{equation*} 
	X_1(t) = \int_0^t \mathcal V(\eta_r) \ud r + \sqrt{2}B_1(t) \; , \; t\ge0 \; \; ,
\end{equation*}
where $\mathcal V(\omega) := \omega^1(0)$, for $\omega = (\omega^1,\omega^2)$. 
Using the so-called \text{Yaglom-reversibility} (see Section 1.4 of \cite{Tothetal2012}), we get that, for every $0 \le s < t$, the random variables $B(t) - B(s)$ and $\int_s^t\mathcal V(\eta_r) \ud r$ are uncorrelated, so that
\begin{equation} \label{set_main_res_yaglom_reversibility}
	\E[X_1(t)^2] = 2 \E|B_1(t)^2] + \E\left[\left(\int_0^t \mathcal V(\eta_r) \ud r\right)^2\right] 
	\; \; .
\end{equation}
This in turn implies that we can rewrite \eqref{set_main_res_def_D} as $D_T(\lambda) = D_B(\lambda) + D_\V(\lambda)$, where for all $\lambda > 0$,
\begin{equation} \label{set_main_res_yaglom_rev_formula}
	D_B(\lambda) := 4\int_0^\infty e^{-\lambda t} \E|B_1(t)^2] \ud t = \dfrac{4}{\lambda^2} \text{ and } D_\V(\lambda) := 2\int_0^\infty e^{-\lambda t} \E\left[\left(\int_0^t \mathcal V(\eta_r) \ud r\right)^2\right] \ud t \; \; ,
\end{equation}
and therefore, we may focus on $D_\V(\lambda)$, which requires a good understanding of the process $(\eta_t)_{t\ge0}$. Since the drift field is stationary (see Proposition \ref{op_fock_prop_invariant_measure_fSHE}) and divergence-free, the law of $\cGFF$ is invariant also for $(\eta_t)_{t\ge0}$ (see e.g. Chapter 11 in \cite{landimetal2012}). This ensures that, by Lemma 5.1 in \cite{CES21}, we can write
\begin{equation} \label{set_main_res_laplace_transform_to_resolvent_eq}
	\int_0^\infty e^{-\lambda t} \E\left[\left(\int_0^t \mathcal V(\eta_r) \ud r\right)^2\right] \ud t = \dfrac{2}{\lambda^2} \Eb[\V(\lambda - \L^s)^{-1}\V] \; \; ,
\end{equation}
where $\L^s$ denotes the infinitesimal generator of $(\eta_t)_{t\ge0}$, defined in \eqref{op_fock_generator_eta} below, and with a slight abuse of notation we use $\Eb$ to denote the expectation with respect to the law of $\cGFF$.

\section{Operators on Fock space} \label{section_operators_fock_space}

In order to analyse expression \eqref{set_main_res_laplace_transform_to_resolvent_eq}, we describe the infinitesimal generator of the infinite dimensional Markov process $t \mapsto \omega_t$. With a small abuse of notation, let $\Pb$ denote the law of $\cGFF$ and consider $F \in L^2(\Pb)$ of the form $F(\omega) = f(\omega^{i_1}(x_1), \dots, \omega^{i_n}(x_n))$ for arbitrary points $x_1, \dots, x_n \in \R^2$ and for
 an $f \in C_p^2(\R^n,\R)$, the $C^2$ functions with polynomially growing partial derivatives of order less or equal than $2$.
In this section, to emphasise its dependence in $s \in [0,\infty)$, let us denote by $\L_0^s$ the infinitesimal generator of $(\omega_t)_{t\ge0}$. For every $s \in [0,\infty)$, an application of It\^o's formula gives
\begin{equation} \label{op_fock_L_0_by_ito}
	\begin{split}
		\L_0^s F(\omega)  &= 
		\sum_{k=1}^{n} \partial_{k}f(\omega^{i_1}(x_1), \dots, \omega^{i_n}(x_n)) (-\flap^s) \omega^{i_k}(x_k) \\
		& + \sum_{k,l=1}^{n} \partial_{kl}^2f(\omega^{i_1}(x_1), \dots, \omega^{i_n}(x_n)) (\flap^{s-1}\cpartial_{i_k}\delta_{x_k}, \cpartial_{i_l}\delta_{x_l})_V,
	\end{split}
\end{equation}
where $(\cdot,\cdot)_V$ is given by \eqref{intro_def_mollified_inner_product}, $\partial_kf$ denotes the function $ y = (y_1,\dots,y_n) \mapsto \partial_{y_k}f(y)$ and for every $x \in \R^2$, the expression with $\delta_x$ is well defined by Remark \ref{intro_remark_mollifier_makes_everything_well_defined}. 

Let us introduce the Wiener chaos with the respect to $\Pb$, following the same convention and notation as \cite{Cannizaroetal2021}. Let $x_{1:n} := (x_1, \dots, x_n)$, $\textbf{i} := (i_1,\dots, i_n)$ and $:\cdots:$ denotes the Wick product with respect to $\Pb$. Define $H_0$ as the set of constant random variables and for $n \geq 1$ let $H_n$ be
the set
\begin{equation}
	\left\{\psi_n = \sum_{i_1,\dots, i_n=1}^{2} \int_{\R^{2n}} f_\textbf{i}(x_{1:n}) :\prod_{k=1}^n \omega^{i_k}(x_k):\ud x_{1:n} \right\} \label{op_fock_def_Hn}
\end{equation}
where
the functions $f_{\textbf{i}}$ are symmetric and such that
\begin{equation} \label{op_fock_convention_fourier}
	\widehat{\psi_n}(p_{1:n}) := (-\iota)^n \sum_{i_1,\dots, i_n=1}^{2} \prod_{k=1}^np^\perp_{k,i_k} \hat{f_\textbf{i}}(p_{1:n}) \; \; ,
\end{equation}
satisfies
\begin{equation} \label{op_fock_finite_variance_fourier}
	\Eb[|\psi_n|^2] = \dfrac{n!}{(2\pi)^{2n}} \int_{\R^{2n}} \prod_{k=1}^n \dfrac{\widehat{V}(p_k)}{|p_k|^2} |\widehat{\psi_n}(p_{1:n})|^2 \ud p_{1:n} < \infty \; \; .
\end{equation}
Here, $(p^\perp_{k,1}, p^\perp_{k,2}) := (p_{k,2}, -p_{k,1})$ for $p_k = (p_{k,1}, p_{k,2})$ and $\hat{f_{\textbf{i}}}$ denotes the Fourier transform of $f_{\textbf{i}}$, given by
\begin{equation*}
	\hat{f_{\textbf{i}}}(p_{1:n}) := \int_{\R^{2n}} f_{\textbf{i}}(x_{1:n}) e^{-\iota x_{1:n} \cdot p_{1:n}} \ud x_{1:n} \; \; ,
\end{equation*}
where $x_{1:n} \cdot p_{1:n}$ denotes the canonical inner product in $\R^{2n}$ and $\iota = \sqrt{-1}$.

\begin{remark} \label{op_fock_remark_V_in_H_1}
	Note that since we have the mollification in the noise, the objects \text{$f_{\textbf{i}}$} can be distributions of any negative regularity, such as the delta Dirac distribution. The random variable which we are most interested in here, namely $\V(\omega) = \omega^1(0)$, defined in the previous section, can be seen as 
	\begin{equation*}
		\V(\omega) = \int_{\R^2} \delta_0(x) \omega^1(x)\ud x \in H_1 \; \; .
	\end{equation*}
	Furthermore, $\hat{\V}(p) = p_2$ for $p = (p_1,p_2)$.
\end{remark}

It is well known, see e.g. Nualart  \cite{Nualart2006} or Janson  \cite{janson1997}, that
\begin{equation} \label{op_fock_def_L2_grading}
	L^2(\Pb) = \bigoplus_{n=0}^\infty H_n
\end{equation}
and for $F^i \in L^2(\Pb), \; i=1,2$ given by $F^i = \sum_{n=0}^{\infty} \psi_n^i$, for $\psi_n^i \in H_n$, the expectation $\Eb[F^1F^2]$ can be written as
\begin{equation} \label{op_fock_def_brackets}
	\Eb[F^1F^2] = \sum_{n=0}^\infty \lb\psi_n^1, \psi_n^2\rb
	:= \sum_{n=0}^\infty
	\dfrac{n!}{(2\pi)^{2n}} \int_{\R^{2n}} \prod_{k=1}^n \dfrac{\widehat{V}(p_k)}{|p_k|^2} \widehat{\psi_n^1}(p_{1:n}) \overline{\widehat{\psi_n^2}(p_{1:n})} \ud p_{1:n} \; \; .
\end{equation}

\begin{remark}
	Henceforth we will implicitly identify a random variable $F \in H_n \subset L^2(\Pb)$ of the form \eqref{op_fock_def_Hn} with its kernel $\widehat{\psi_n}$ in Fourier space.  In the same philosophy, we will denote linear operators acting on $L^2(\Pb)$ with the correspondent operators acting on Fock space $\bigoplus_n L^2_{sym}(\R^{2n})$, and we will denote them by the same symbol.
\end{remark}

Now we are ready to prove

\begin{proposition} \label{op_fock_prop_invariant_measure_fSHE}
	The action of the infinitesimal generator $\L_0^s$ in \eqref{op_fock_L_0_by_ito} is diagonal in Fock space ($\L_0^s : H_n \to H_n$), and is given by
	\begin{equation} \label{op_fock_action_L_o_wick_monomials}
		\L_0^s :\omega^{i_1}(x_1) \cdots \omega^{i_n}(x_n):  \; \; 
		= \sum_{k=1}^{n} :\omega^{i_1}(x_1) \cdots (-\flap^s)\omega^{i_k}(x_k) \cdots \omega^{i_n}(x_n):
	\end{equation}
	on Wick monomials, and in Fourier variables by
	\begin{equation} \label{op_fock_action_L_o_fourier}
		\widehat{(-\L_0^s)\psi_n}(p_{1:n}) = \sum_{k=1}^{n} |p_k|^{2s} \widehat{\psi_n}(p_{1:n}) \; \; .
	\end{equation}
	Furthermore, the law of $\cGFF$ is invariant under the dynamics governed by $\L_0^s$, i.e., the infinite dimensional Markov process $(\omega_t)_{t\ge0}$ is stationary and it is distributed according to the law of $\cGFF$ for every $t \ge 0$.
\end{proposition}
\begin{proof}
	By the definition of Wick monomials, we have that 
	\begin{equation*}
		\partial_{k}:\omega^{i_1}(x_1) \cdots \omega^{i_n}(x_n): \;  = \; :\omega^{i_1}(x_1) \cdots \cancel{\omega^{i_k}(x_k)} \cdots \omega^{i_n}(x_n): \; \; ,
	\end{equation*}
	where $a\cancel{b}c := ac$ for $a,b,c \in \R$. Now, the above applied to \eqref{op_fock_L_0_by_ito} with $F = \; :\omega^{i_1}(x_1) \cdots \omega^{i_n}(x_n):$ gives
	\begin{align}
			\L_0^s F(\omega)  &= 
			\sum_{k=1}^{n} :\omega^{i_1}(x_1) \cdots \cancel{\omega^{i_k}(x_k)} \cdots \omega^{i_n}(x_n): (-\flap^s) \omega^{i_k}(x_k) \label{op_fock_prop_inv_mea_proof_mult_by_w}\\
			& + \sum_{k,l=1 \atop k \neq l}^{n} :\omega^{i_1}(x_1) \cdots \cancel{\omega^{i_k}(x_k)} \cdots \cancel{\omega^{i_l}(x_l)} \cdots \omega^{i_n}(x_n): (\flap^{s-1}\cpartial_{i_k}\delta_{x_k}, \cpartial_{i_l}\delta_{x_l})_V \; \; . \label{op_fock_prop_inv_mea_proof_lower_chaos_term}
	\end{align}
	Note that on Wick monomials, multiplication by $\omega^{i_k}(x_k)$, as in \eqref{op_fock_prop_inv_mea_proof_mult_by_w}, produces both a term in one higher homogeneous chaos and a term in one lower homogeneous chaos. Precisely, for each $1 \le k \le n$ in \eqref{op_fock_prop_inv_mea_proof_mult_by_w} we have
	\begin{align*}
		& :\omega^{i_1}(x_1) \cdots \cancel{\omega^{i_k}(x_k)} \cdots \omega^{i_n}(x_n): (-\flap^s) \omega^{i_k}(x_k) = \\
		& :\omega^{i_1}(x_1) \cdots (-\flap^s)\omega^{i_k}(x_k) \cdots \omega^{i_n}(x_n): \\
		& + \sum_{l=1 \atop l \neq k}^{n} :\omega^{i_1}(x_1) \cdots \cancel{\omega^{i_k}(x_k)} \cdots \cancel{\omega^{i_l}(x_l)} \cdots \omega^{i_n}(x_n): (\flap^{-1}(-\flap^s)\cpartial_{i_k}\delta_{x_k}, \cpartial_{i_l}\delta_{x_l})_V \; \; ,
	\end{align*}
	where $(\flap^{-1}(-\flap^s)\cpartial_{i_k}\delta_{x_k}, \cpartial_{i_l}\delta_{x_l})_V = \Eb[(-\flap^s) \omega^{i_k}(x_k), \omega^{i_l}(x_l)]$. Summing over $k$, the first term after the equal sign gives us \eqref{op_fock_action_L_o_wick_monomials} and the second term after the equality cancels out with \eqref{op_fock_prop_inv_mea_proof_lower_chaos_term}. \eqref{op_fock_action_L_o_fourier} is a direct consequence of \eqref{op_fock_action_L_o_wick_monomials} and \eqref{op_fock_convention_fourier}. Now we move to the invariance of the law of $\cGFF$. It is known that a necessary and sufficient condition for this is that $(\L_0^{s})^* \textbf{1} = 0$, where $(\L_0^{s})^*$ denotes the adjoint of the operator $\L_0^s$ in $L^2(\Pb)$
	and $\textbf{1}$ denotes the constant function equal to $1$, see e.g. \cite[Theorem 3.37]{Liggett2010}. Also, by \eqref{op_fock_def_L2_grading}, it is enough to consider $F = \; :\omega^{i_1}(x_1) \cdots \omega^{i_n}(x_n):$ , so that
	\begin{equation*}
		\Eb[(\L_0^{s})^* \textbf{1} F] = \Eb[\L_0^sF] = \sum_{k=1}^{n} \Eb[:\omega^{i_1}(x_1) \cdots (-\flap^s)\omega^{i_k}(x_k) \cdots \omega^{i_n}(x_n):] = 0
	\end{equation*}
	completes the proof.
\end{proof}

So far we gathered all the ingredients necessary to characterise the full generator $\L =: \L^s$ of $(\eta_t)_{t\ge0}$. Putting together the generator $\L_0^s$ of the environmental process $(\omega_t)_{t\ge0}$ with Proposition \ref{op_fock_prop_invariant_measure_fSHE}, the arguments in Section 2.1 of T{\'o}th and Valk{\'o} \cite{TothandValko2012} and the main result of Komorowski \cite{Komorowski2000}, we get that the generator $\L^s$ is given by
\begin{equation} \label{op_fock_generator_eta}
	\L^s = \L_0^s + \mathcal{V}\nabla + \Delta = \L_0^s + \A+ - \A+^* + \Delta \; \; ,
\end{equation}
where $\mathcal{V}\nabla := \mathcal{V}_1D_1 + \mathcal{V}_2D_2$, with $\mathcal{V}_i(\omega) = \omega^i(0)$ and $D_i$ is the infinitesimal generator of the spatial shifts in the canonical directions of $\R^2$, for $i=1,2$, see \cite{Komorowski2000}. Also, $\mathcal{V}\nabla = \A+ - \A+^*$ can be decomposed into a creation and annihilation parts, one being minus the adjoint of the other, and it comes from the drift part of \eqref{intro_def_main_SDE}, i.e., the environment, while $\Delta = \nabla^2$ comes from the Brownian part in \eqref{intro_def_main_SDE}, see \cite{TothandValko2012}. We have that
\begin{equation*}
	\L_0^s, \Delta : H_n \to H_n \quad, \quad \A+ : H_n \to H_{n+1} \quad \text{and} \quad \A+^* : H_n \to H_{n-1} \; \; .
\end{equation*}


As noted in Toninelli et al.  \cite{Cannizaroetal2021}, adopting the conventions on Fock space discussed earlier, one has
\begin{align}
	& \widehat{(-\Delta)\psi_n}(p_{1:n}) = \left|\sum_{k=1}^{n} p_k\right|^{2} \widehat{\psi_n}(p_{1:n}) \; \; , \label{op_fock_action_Delta} \\
	& \widehat{\A+\psi_n}(p_{1:n+1}) = \iota \dfrac{1}{n+1} \sum_{k=1}^{n+1} \left(p_k \times \sum_{l=1}^{n+1}p_l\right)\widehat{\psi_n}(p_{1:n+1\setminus l}) \; \; ,
\end{align}
where $p_{1:n+1\setminus l} := (p_1, \dots, \cancel{p_l}, \dots, p_{n+1})$ and for $p,q \in \R^2$, $p \times q$ denotes the scalar given by the third coordinate of the cross product of $p$ with $q$, when thought as vectors in $\R^3$, precisely, $p \times q = p_1q_2 - p_2q_1 = |p||q| \sin \theta$, where $\theta$ is the angle between $p$ and $q$.

\begin{remark} \label{op_fock_remark_heuristics_for_main_theo}
	Here we can see that if $s=1$ in \eqref{op_fock_action_L_o_fourier}, the difference between the operators $\Delta$ and $\L_0^1$ is simply the cross terms in \eqref{op_fock_action_Delta}. The most important observation here is that if $s\ge1$ and $|p|\le1$, in view of \eqref{set_main_res_laplace_transform_to_resolvent_eq} and Remark \ref{op_fock_remark_V_in_H_1},
	for any function $\psi_1 \in H_1$, we have that
	\begin{equation*}
		\widehat{(-\Delta)\psi_1}(p) = |p|^{2}\widehat{\psi_1}(p) \ge |p|^{2s} \widehat{\psi_1}(p) = \widehat{(-\L_0^s)\psi_1}(p) \; \; .
	\end{equation*}
	This is a good evidence to suggest \eqref{intro_super_dif_bounds_s>1},
	 as can be further seen in Remark \ref{anal_res_remark_level_1}. Also, a good heuristics for the drastic change in behaviour in $s$ contained in Theorem \ref{intro_theo_main_theo} is that in Fourier variables, the operator $\L_0^s$ acts much more severely in large scales when $s<1$ than when $s\ge1$, since $|p|^{2s} << |p|^{2s^\prime}$ for $|p| << 1$, if $s^\prime < s$.
\end{remark}

Now we proceed to the analysis of the resolvent equation in \eqref{set_main_res_laplace_transform_to_resolvent_eq}.

\section{Iterative analysis of the resolvent equation} \label{section_iterative_anal_resolvent_eq}

We can write $\Eb[\V(\lambda - \L^s)^{-1}\V]$ as $\Eb[\V\mathbf V]$, where $\mathbf V$ is the solution to the resolvent equation $(\lambda - \L^s)\mathbf V = \V$. Note however that $\V \in H_1$ is in the first Wiener chaos and that the operator $\L^s$ maps $H_n$ to $H_{n-1} \oplus H_{n} \oplus H_{n+1}$, one should expect that the solution $\mathbf V$ to the resolvent equation has non-trivial componentes in all Wiener chaoses.
Following the idea introduced by Landim et al. \cite{LQSY2004} we truncate the generator $\L^s$ by using $\L_n^s := P_{\le n} \L^s P_{\le n}$, where $P_{\le n}$ denotes the orthogonal projection onto the inhomogeneous chaos of order $n$, i.e., $P_{\le n} : L^2(\Pb) \to \bigoplus_{k=0}^n H_k$. Denote by $\Vs^n \in \bigoplus_{k=0}^n H_k$ the solution to the resolvent equation truncated at level $n$, i.e.,
\begin{equation*}
	(\lambda - \L_n^s)\Vs^n = \V \text{ and } \Vs^n = \sum_{k=0}^n \Vs^n_k \text{ where } \Vs^n_k \in H_k \; , \; k = 0,1,\dots, n \; \; .
\end{equation*}
Now, writing one equation for each of the components of $\Vs$ above we get that the equation above is equivalent to the system of equations
\begin{equation*}
	\begin{cases}
		& (\lambda - \Delta - \L_0^s)\Vs^n_n - \A+\Vs^n_{n-1} = 0 \; , \\
		& (\lambda - \Delta - \L_0^s)\Vs^n_{n-1} - \A+\Vs^n_{n-2} + \A+^*\Vs^n_{n} = 0 \; , \\
		& \cdots \\
		& (\lambda - \Delta - \L_0^s)\Vs^n_{1} + \A+^*\Vs^n_{2} = \V \; \; ,
	\end{cases}
\end{equation*}
Note that as it was observed in \cite[Section 2]{TothandValko2012}, $\A+^*F = 0$ for every $F \in H_1$, so that $\Vs^n_0 = 0$ and we do not write an equation for it. Note that since $\V \in H_1$ to evaluate \eqref{set_main_res_laplace_transform_to_resolvent_eq} at the level of the truncation, only the component in the first Wiener chaos is necessary, i.e., $\Vs^n_1$. For that, the system above can be solved and shows that 
\begin{equation*}
	\Vs^n_1 = (\lambda - \Delta - \L_0^s + \H_n)^{-1}\V \; \; ,
\end{equation*}
where
\begin{equation} \label{anal_res_eq_def_H_k}
	\begin{cases}
		& \H_1 := 0 \; , \\
		& \H_{k+1} = \A+^*(\lambda - \Delta - \L_0^s + \H_k)^{-1} \A+ \; \; , k \ge 1 \; \; .
	\end{cases}
\end{equation}
It is important to note that $\H_k : H_n \to H_n$ for every $k,n \in \N$. Recall that by \eqref{op_fock_def_brackets} we can write $\Eb[\V(\lambda - \L_n^s)^{-1}\V] = \lb \V , \Vs^n_1 \rb$. As it was first noticed in \cite[eq. (2.4)]{LQSY2004}, the following monotonicity formula follows from the fact that $\lambda - \Delta - \L_0^s$ is a positive operator.
\begin{lemma} \label{anal_res_eq_even_odd_inequalities}
	Let $S := \lambda - \Delta - \L_0^s$, then, for every $n \ge 1$, we get the bounds
	\begin{equation*}
		\lb \V , (S + \H_{2n})^{-1}\V \rb = \lb \V , \Vs^{2n}_1 \rb \le \lb \V , (\lambda - \L^s)^{-1}\V \rb \le \lb \V , \Vs^{2n-1}_1 \rb = \lb \V , (S + \H_{2n-1})^{-1}\V \rb \; \; .
	\end{equation*}
\end{lemma} 

\begin{remark} \label{anal_res_remark_level_1}
	Let us look to the first upper bound when taking $n=1$ in Lemma \ref{anal_res_eq_even_odd_inequalities} above. Recall that $\V \in H_1$ and that $\hat{\V}(p) = p_2$ for $p = (p_1,p_2)$. Thus by considering the solution $\Vs^1$ to the truncation at the first level, we arrive at
	\begin{align}
		\lb \V , (\lambda - \L^s)^{-1}\V \rb 
		\le \lb \V , (\lambda - \Delta - \L_0^s)^{-1}\V \rb
		& = \frac{1}{(2\pi)^2} \int_{\R^2} \dfrac{\hat{V}(p)}{|p|^2} \dfrac{|\hat{\V}(p)|^2 \ud p}{\lambda + |p|^2 + |p|^{2s}} \nonumber\\
		& \le C \int_0^1\dfrac{r \ud r}{\lambda + r^2 + r^{2s}}
		\le
		\begin{cases}
			& C \quad \quad \quad \quad \quad \quad \text{if } s<1 \; , \\
			& C \log\left(1 + \dfrac{1}{\lambda}\right) \text{ if } s \ge 1\; ,
		\end{cases}
		\label{anal_res_remark_level_1_eq}
	\end{align}
	for a constant $C>0$. Note now that for the case $s<1$, the inequalities in \eqref{anal_res_remark_level_1} imply the diffusive bounds \eqref{intro_dif_bounds_s<1} in Theorem \ref{intro_theo_invariance_principle}, see \eqref{proof_theo_inv_princ_H_-1_norm_bound_variance} in Section \ref{section_proof_theo_inv_princ} and the following discussion. On the other hand, for the case of $s \ge 1$, the estimates in \eqref{anal_res_remark_level_1} together with the first lower bound obtained with $n=1$ in Lemma \ref{anal_res_eq_even_odd_inequalities}, by the same argument for the lower bound in Section \ref{section_proof_theo_log_div} for the case $\gamma=1$, gives
	\begin{equation} \label{intro_non_sharp_estimates_TV12}
		A \log |\log \lambda| \le \lambda^2 D_T(\lambda) \le B |\log \lambda| \text{ for } \lambda \in (0,1) \; \; ,
	\end{equation}
	for constants $A,B>0$. These are precisely the estimates obtained in \cite{TothandValko2012} for the static case. In particular this already implies that the dynamics of SHE is not enough to remove the super diffusivity caused by the random environment.
\end{remark}

The estimates 
in \eqref{intro_non_sharp_estimates_TV12} 
can be iterated for higher levels and be improved at each step. Indeed, to get 
\eqref{intro_super_dif_bounds_s>1},
it is necessary to use Lemma \ref{anal_res_eq_even_odd_inequalities} in full by taking the level $k$ to diverge with $\lambda \downarrow 0$. Moreover, an understanding of the estimates for every level is necessary, and for that
it suffices to analyse the operators $\H_k$. For this, we make use of the following three lemmas, taken from Toninelli et al. \cite{Cannizaroetal2021}. In what follows, $\mathcal S$ is an operator which acts diagonally in Fock space with Fourier multiplier denoted by $\sigma$, such that $\widehat{\S\psi_n}(p_{1:n}) = \sigma_n(p_{1:n})\widehat{\psi_n}(p_{1:n})$ for any $\psi_n \in H_n$, which will later be taken to be $\mathcal S = S + \H_n$, for $n \ge 1$.

\begin{lemma}\label{anal_res_CHT_lemma_44} 
	For any $\psi_n \in H_n$, it holds that
	\begin{equation*}
		\lb \psi_n, \A+^*\S\A+\psi_n \rb = \lb \psi_n, \A+^*\S\A+\psi_n \rb_{\text{Diag}} + \lb \psi_n, \A+^*\S\A+\psi_n \rb_{\text{Off}} \; \; ,
	\end{equation*}
	where
	\begin{equation*}
		\lb \psi_n, \A+^*\S\A+\psi_n \rb_{\text{Diag}}
		:= \dfrac{n!}{(2\pi)^{n+1}} \int_{\R^{2(n+1)}} \prod_{k=1}^{n+1} \dfrac{\widehat{V}(p_k)}{|p_k|^2} |\widehat{\psi_n}(p_{1:n})|^2 \sigma_{n+1}(p_{1:n+1}) \left(p_{n+1} \times \sum_{k=1}^{n}p_k\right)^2 \ud p_{1:n+1}
	\end{equation*}
	and
	\begin{align*}
		&\lb \psi_n, \A+^*\S\A+\psi_n \rb_{\text{Off}} \\
		&:= \dfrac{n!n}{(2\pi)^{n+1}} \int_{\R^{2(n+1)}} \prod_{k=1}^{n+1} \dfrac{\widehat{V}(p_k)}{|p_k|^2} \overline{\widehat{\psi_n}(p_{1:n})}\widehat{\psi_n}(p_{1:n+1\setminus n}) \sigma_{n+1}(p_{1:n+1}) \left(p_{n+1} \times \sum_{k=1}^{n+1}p_k\right)\left(p_{n} \times \sum_{k=1}^{n+1}p_k\right) \ud p_{1:n+1}
	\end{align*}
\end{lemma}

\begin{lemma} \label{anal_res_CHT_lemma_45}
	If for every $n \in \N$ and any $p_{1:n} \in \R^{2n}$ with $\sum_{k=1}^{n} p_k \neq 0$
	\begin{equation} \label{anal_res_lemmas_CHT_Diag_1}
		\int_{\R^2} \widehat{V}(q) (\sin \theta)^2 \sigma_{n+1}(p_{1:n},q) \ud q \le \tilde{\sigma}_n(p_{1:n})
	\end{equation}
	with $\theta$ the angle between $q$ and $\sum_{k=1}^{n} p_k$, then for every $\psi_n$
	\begin{equation} \label{anal_res_lemmas_CHT_Diag_2}
		\lb \psi_n, \A+^*\S\A+\psi_n \rb_{\text{Diag}} \le \lb \psi_n, (-\Delta)\tilde{\S}\psi_n \rb
	\end{equation}
	where $\tilde{\S}$ is the diagonal operador whose Fourier multiplier is $\tilde{\sigma}$. If the inequality in \eqref{anal_res_lemmas_CHT_Diag_1} is $\ge$, then \eqref{anal_res_lemmas_CHT_Diag_2} holds with $\ge$ as well.
\end{lemma}

\begin{lemma} \label{anal_res_CHT_lemma_46}
	If for every $n \in \N$ and any $p_{1:n} \in \R^{2n}$
	\begin{equation*}
		\left|\sum_{k=1}^{n} p_k\right|\int_{\R^2} \widehat{V}(q) \dfrac{(\sin \theta)^2 \sigma_{n+1}(p_{1:n},q)}{\left|q + \sum_{k=1}^{n-1} p_k\right|} \ud q \le \tilde{\sigma}_n(p_{1:n})
	\end{equation*}
	with $\theta$ the angle between $q$ and $\sum_{k=1}^{n} p_k$, then for every $\psi_n$
	\begin{equation*}
		|\lb \psi_n, \A+^*\S\A+\psi_n \rb_{\text{Off}}| \le n \lb \psi_n, (-\Delta)\tilde{\S}\psi_n \rb
	\end{equation*}
	where $\tilde{\S}$ is the diagonal operador whose Fourier multiplier is $\tilde{\sigma}$.
\end{lemma}

Here are some preliminary definitions, needed to state and prove the next theorem. Expressions \eqref{anal_fock_def_L} and \eqref{anal_fock_def_LB_UB} arise naturally when iterating the estimates for different levels $k$ in Lemma \ref{anal_res_eq_even_odd_inequalities}.
For $k \in \N, x > 0$ and $z \ge 0$, let $L, LB_k$ and $UB_k$ be given by
\begin{equation} \label{anal_fock_def_L}
	L(x,z) := z + \log(1 + x^{-1}) \; ,
\end{equation}

\begin{equation} \label{anal_fock_def_LB_UB}
	LB_k(x,z) := \sum_{j=0}^{k} \dfrac{(1/2 \log L(x,z))^j}{j!} \text{  and  } UB_k(x,z) := \dfrac{L(x,z)}{LB_k(x,z)}
\end{equation}

\noindent and, for $k \ge 1$, define $\sigma_k$ as

\begin{equation*}
	\sigma_k(x,z) = 
	\begin{cases}
		& UB_{\frac{k-2}{2}}(x,z) \,, \text{ if } k \text{ is even,} \; \\
		& LB_{\frac{k-1}{2}}(x,z) \,, \text{ if } k \text{ is odd.}
	\end{cases}
\end{equation*}

We have that $\sigma_k \equiv 1$. Also, for $n \in \N$, let

\begin{equation} \label{anal_fock_def_zk_fk}
	z_k(n) = K_1(n + k)^{2 + 2\varepsilon} \text{  and  } f_k(n) = K_2\sqrt{z_k(n)} \; \; ,
\end{equation}
where $K_1,K_2$ are constants to be chosen sufficiently large later and $\varepsilon$ is the small positive constant appearing in the main Theorem \ref{intro_theo_main_theo}. Now, for $k \ge 1$, let $\delta_k$ be an operator such that its Fourier multiplier is $\sigma_k$, meaning

\begin{equation*}
	\delta_k = 
	\begin{cases}
		& f_k(\Nop)\sigma_k(\lambda - \Delta - \L_0^s, z_k(\Nop)) \,, \text{ if } k \text{ is even,} \; \\
		& \frac{1}{f_k(\Nop)}(\sigma_k(\lambda - \Delta - \L_0^s, z_k(\Nop)) - f_k(\Nop)) \,, \text{ if } k \text{ is odd.}
	\end{cases}
\end{equation*}
where $\Nop$ denotes the so-called Number Operator, the infinitesimal generator of $\partial_tu = -u + \sqrt{2}(-\Delta)^{-\frac{1}{2}}\xi$, which acts diagonally on the $n$-th Wiener chaos by multiplying by $n$: $\Nop \psi_n = n \psi_n$ for every $\psi_n \in H_n$.

\begin{remark}
	Note that the functions $L, LB_k$ and $UB_k$ are the same as in Toninelli et al. \cite{Cannizaroetal2021}, while the operators $\delta_k$ carry the the generator $\L_0^s$, which is the difference between the dynamic and the static settings.
\end{remark}

Gathering these we put them into the next theorem.

\begin{theorem} \label{anal_fock_theo_op_estimates}
	If $s\ge1$ in \eqref{intro_def_environmental_process_omega}, then for every $\varepsilon > 0$, we may choose $K_1$ and $K_2$ in \eqref{anal_fock_def_zk_fk} to be large enough so that, for $0 < \lambda \le 1$ and $k \ge 1$, the following operator estimates hold true.
	\begin{equation} \label{anal_fock_theo_op_estimates_lower_s>1}
		\H_{2k-1} \ge c_{2k-1} (-\Delta)\delta_{2k-1}
	\end{equation}
	and
	\begin{equation} \label{anal_fock_theo_op_estimates_upper_s>1}
		\H_{2k} \le c_{2k} (-\Delta)\delta_{2k} \; .
	\end{equation}
	where $c_1 = 1$ and
	\begin{equation} \label{anal_fock_theo_op_estimates_ck}
		c_{2k} = \dfrac{\pi}{c_{2k-1}}\left(1 + \dfrac{1}{k^{1 + \varepsilon}}\right) \;, \; \; 
		c_{2k+1} = \dfrac{\pi}{c_{2k}}\left(1 - \dfrac{1}{(k+1)^{1 + \varepsilon}}\right) \; .
	\end{equation}	
\end{theorem}

\begin{remark}
	We shall emphasise here that the sequences $c_{2k}$ and $c_{2k+1}$ in \eqref{anal_fock_theo_op_estimates_ck} do converge to finite, strictly positive constants, as $k \to \infty$, provided that $\varepsilon > 0$. Furthermore, the limits are strictly greater than $2\pi$ and strictly smaller than $1$, respectively. This can be seen, e.g. for the even sequence, $\frac{c_{2k+2}}{c_{2k}} = (1 + \frac{1}{k^{1 + \varepsilon}})(1 - \frac{1}{(k+1)^{1 + \varepsilon}})^{-1} > 1$ and $c_2 = 2\pi$. Also, by iterating the definition for $c_{2k}$, it can be shown that convergence of the sequence is equivalent to the convergence of $\sum_{l=1}^\infty l^{-(1+\varepsilon)}$, which only holds when $\varepsilon>0$.
\end{remark}

Now we will prove Theorem \ref{anal_fock_theo_op_estimates} by induction on $k$.
Note that the induction alternates between lower \eqref{anal_fock_theo_op_estimates_lower_s>1} and upper \eqref{anal_fock_theo_op_estimates_upper_s>1} bounds, being one the consequence of the other, and so forth. 

\begin{proof}[Proof of the lower bound \eqref{anal_fock_theo_op_estimates_lower_s>1}]
	Recall that $s \ge 1$. For $k=1$ we note that, by definition, $\H_1 = 0$ and $\delta_1$ is non-positive if we choose the constant $K_2$ in \eqref{anal_fock_def_zk_fk} to be large enough. 
	
	We now show \eqref{anal_fock_theo_op_estimates_lower_s>1} with $2k+1$ for $k\ge1$, assuming by induction that \eqref{anal_fock_theo_op_estimates_upper_s>1} holds for $2k$:
	\begin{equation} \label{anal_fock_proof_op_est_lower_start_H_2k+1}
		\H_{2k+1} = \A+^*(\lambda - \Delta - \L_0^s + \H_{2k})^{-1} \A+ \ge \A+^*(\lambda - \Delta(1 + c_{2k}\delta_{2k}) - \L_0^s)^{-1} \A+ \; .
	\end{equation}
	For every $\psi \in H_n$, we use Lemma \ref{anal_res_CHT_lemma_44} with $S = (\lambda - \Delta(1 + c_{2k}\delta_{2k}) - \L_0^s)^{-1}$ to separate
	\begin{equation} \label{anal_fock_proof_op_est_lower_diag_+_off}
		\lb\psi, \A+^*(\lambda - \Delta(1 + c_{2k}\delta_{2k}) - \L_0^s)^{-1} \A+ \psi\rb = \lb\psi, \A+^*S \A+ \psi\rb
	\end{equation}
	into a diagonal and an off-diagonal part, and we treat each separately. For the diagonal part, we apply Lemma \ref{anal_res_CHT_lemma_45} for which it suffices to lower bound
	\begin{equation} \label{anal_fock_proof_op_est_lower_int_1}
		\int_{\R^2} \dfrac{\hat{V}(q)(\sin \theta)^2 \ud q}{\lambda + |p+q|^2(1 + c_{2k}f_{2k}(n+1)UB_{k-1}(\lambda + |p+q|^2 + |p_{1:n}|^{2s} + |q|^{2s},z_{2k}(n+1))) + |p_{1:n}|^{2s} + |q|^{2s}}
	\end{equation}
	where $p = \sum_{i=1}^{n} p_i$ and $|p_{1:n}|^{2s} := \sum_{i=1}^{n} |p_i|^{2s}$, for $p_1,\dots, p_n \in \R^2$, and $\theta$ is the angle between $p$ and $q$. Clearly, $|p_{1:n}|^{2s}$ is different from $|p|^2$ even for $s=1$. Naturally, the argument in $z_{2k},f_{2k}$ is $n+1$ since $\A+\psi \in H_{n+1}$, but by \eqref{anal_fock_def_zk_fk} we get that $z_{2k}(n+1) = z_{2k+1}(n)$ and $f_{2k}(n+1) = f_{2k+1}(n)$ and henceforth we drop the argument $n$ to lighten the notation. We may upper bound the denominator in \eqref{anal_fock_proof_op_est_lower_int_1} by
	\begin{equation*}
		(\lambda + |p+q|^2 + |p_{1:n}|^{2s})(1 + c_{2k}f_{2k+1}UB_{k-1}(\lambda + |p+q|^2 + |p_{1:n}|^{2s} + |q|^{2s},z_{2k+1}))  + |q|^{2s}
	\end{equation*}
	\begin{equation*}
		\le c_{2k}f_{2k+1}\left(1 + \frac{1}{f_{2k+1}}\right)[(\lambda + |p+q|^2 + |p_{1:n}|^{2s})UB_{k-1}(\lambda + |p+q|^2 + |p_{1:n}|^{2s},z_{2k+1})  + |q|^{2s}] \; ,
	\end{equation*}
	where we have used for both inequalities that $c_{2k}, f_{2k+1}, UB_{k-1} \ge 1$ and the monotonicity of $UB_{k-1}$. Therefore, we may look to
	\begin{equation*}
		\int_{\R^2} \dfrac{\hat{V}(q)(\sin \theta)^2 \ud q}{(\lambda + |p_{1:n}|^{2s} + |p+q|^2)UB_{k-1}(\lambda + |p_{1:n}|^{2s} + |p+q|^2,z_{2k+1}) + |q|^{2s}} \; \; .
	\end{equation*}
	By 
	Lemma \ref{appendix_lemma_3}, we get for the integral above, the lower bound
	\begin{align}
		& \dfrac{\pi}{2}\int_{\lambda + |p|^2 + |p_{1:n}|^{2s}}^1 \dfrac{\ud \varrho}{\varrho UB_{k-1}(\varrho,z_{2k+1})} - 
		C_{\text{Diag}} \dfrac{LB_{k}(\lambda + |p|^2 + |p_{1:n}|^{2s},z_{2k+1})}{\sqrt{z_{2k+1}}} \label{anal_fock_proof_op_est_lower_dominant_with_pi/4} \\
		& \ge \dfrac{\pi}{2}\int_{\lambda + |p|^2 + |p_{1:n}|^{2s}}^1 \dfrac{\ud \varrho}{(\varrho + \varrho^2) UB_{k-1}(\varrho,z_{2k+1})} - 
		C_{\text{Diag}} \dfrac{LB_{k}(\lambda + |p|^2 + |p_{1:n}|^{2s},z_{2k+1})}{\sqrt{z_{2k+1}}} \; . \nonumber
	\end{align}
	By \eqref{appendix_CHT_exp_A2} the primitive of the integral above is $-2LB_k(\varrho,z_{2k+1})$, hence the expression above equals
	\begin{align*}
		& \pi LB_k(\lambda + |p|^2 + |p_{1:n}|^{2s},z_{2k+1}) - \pi LB_k(1,z_{2k+1}) - 
		C_{\text{Diag}} \dfrac{LB_{k}(\lambda + |p|^2 + |p_{1:n}|^{2s},z_{2k+1})}{\sqrt{z_{2k+1}}} \\
		& \ge \pi LB_k(\lambda + |p|^2 + |p_{1:n}|^{2s},z_{2k+1}) - \dfrac{\pi f_{2k+1}}{2} - 
		C_{\text{Diag}} \dfrac{LB_{k}(\lambda + |p|^2 + |p_{1:n}|^{2s},z_{2k+1})}{\sqrt{z_{2k+1}}}
	\end{align*}
	where in the last inequality we have again used Lemma \ref{appendix_CHT_lemma_A1} and chosen the constant $K_2$ in \eqref{anal_fock_def_zk_fk} large enough so that for all $k,n \in \N$, it holds that
	\begin{equation} \label{anal_fock_proof_op_est_lower_ineq_LB(1,z)}
		LB_k(1,z_{2k+1}) \le \sqrt{L(1,z_{2k+1})} = \sqrt{\log(2) + z_{2k+1}} \le \frac{1}{2}f_{2k+1} \; \; .
	\end{equation}
	So by Lemma \ref{anal_res_CHT_lemma_45} we get that the diagonal part of \eqref{anal_fock_proof_op_est_lower_diag_+_off} is lower bounded by $\lb\psi, (-\Delta) \tilde{S} \psi\rb$, where
	\begin{equation} \label{anal_fock_proof_op_est_lower_diag_S_tilde}
		\tilde{S} = \left(1 + \dfrac{1}{f_{2k+1}(1)}\right)^{-1} \dfrac{\pi}{c_{2k}}
		\left[\dfrac{LB_k(\lambda - \Delta - \L_0^s, z_{2k+1}(\Nop))}{f_{2k+1}(\Nop)}\left(1 - \dfrac{C_{\text{Diag}}}{\pi \sqrt{z_{2k+1}(1)}}\right) - \dfrac{1}{2}\right] \; \; .
	\end{equation}
	Here we have twice lower bounded $z_{2k+1} = z_{2k+1}(n) \ge z_{2k+1}(1)$ and $f_{2k+1} = f_{2k+1}(n) \ge f_{2k+1}(1)$.
	
	For the off-diagonal part of \eqref{anal_fock_proof_op_est_lower_diag_+_off} we use Lemma \ref{anal_res_CHT_lemma_46}. For that, denote $p = \sum_{i=1}^{n} p_i$ and $p^\prime = \sum_{i=1}^{n-1} p_i$ and we must upper bound
	\begin{align} \label{anal_fock_proof_op_est_lower_int_2}
		& n|p| \int_{\R^2} \dfrac{\hat{V}(q)(\sin \theta)^2 \ud q}{[\lambda + |p+q|^2(1 + c_{2k}f_{2k+1}UB_{k-1}(\lambda + |p+q|^2 + |p_{1:n}|^{2s} + |q|^{2s},z_{2k+1})) + |p_{1:n}|^{2s} + |q|^{2s}] |p^\prime + q|} \nonumber \\
		& \le n|p| \int_{\R^2} \dfrac{\hat{V}(q)(\sin \theta)^2 \ud q}{[\lambda + |p_{1:n}|^{2s} + |p+q|^2c_{2k}f_{2k+1}UB_{k-1}(\lambda + |p_{1:n}|^{2s} + |p+q|^2 + |q|^{2s},z_{2k+1})] |p^\prime + q|} \nonumber \\
		& \le n|p| \int_{\R^2} \dfrac{\hat{V}(q)(\sin \theta)^2 \ud q}{[\lambda + |p_{1:n}|^{2s} + |p+q|^2c_{2k}f_{2k+1}UB_{k-1}(\lambda + |p_{1:n}|^{2s} + |p+q|^2 + |q|^{2},z_{2k+1})] |p^\prime + q|} \; \; ,
	\end{align}
	where in the last inequality we have used the monotonicity of $UB_{k-1}$ and that since $\hat{V}$ is supported on $|q| \le 1$, we have that $|q|^{2s} \le |q|^2$ if $s\ge1$.Thanks to Lemma \ref{appendix_CHT_lemma_A4} the functions $f(x,z) = c_{2k}f_{2k+1}UB_{k-1}(x,z)$ and $g(x,z) = \frac{1}{c_{2k}f_{2k+1}}LB_{k-1}(x,z)$ satisfy the assumptions of Lemma \ref{appendix_lemma_2} and we obtain the upper bound
	\begin{equation} \label{anal_fock_proof_op_est_lower_off_S_tilde}
		\dfrac{n C_{\text{Off}} LB_{k-1}(\lambda + |p_{1:n}|^{2s} + |p|^2,z_{2k+1})}{c_{2k}f_{2k+1}z_{2k+1}}
		\le \dfrac{C_{\text{Off}}}{c_{2k}f_{2k+1}K_1(2k + 1)^{1+\varepsilon}}  LB_{k}(\lambda + |p_{1:n}|^{2s} + |p|^2,z_{2k+1}) 		
	\end{equation}
	where we have used that $LB_{k-1} \le LB_{k} $, the definition of $z_{2k+1} = z_{2k+1}(n)$ in \eqref{anal_fock_def_zk_fk} and the fact that
	\begin{equation*}
		\dfrac{n}{z_{2k+1}(n)} = \dfrac{n}{K_1(2k + 1 + n)^{2+2\varepsilon}} \le \dfrac{1}{K_1(2k + 1 + n)^{1+\varepsilon}} \; \; .
	\end{equation*}
	
	Altogether, Lemmas \ref{anal_res_CHT_lemma_45} and \ref{anal_res_CHT_lemma_46} combined with expressions \eqref{anal_fock_proof_op_est_lower_diag_S_tilde} and \eqref{anal_fock_proof_op_est_lower_off_S_tilde}, we obtain that the operator $\A+^*(\lambda - \Delta(1 + c_{2k}\delta_{2k}) - \L_0^s)^{-1} \A+$ is lower bounded by
	\begin{equation*}
		(-\Delta)\dfrac{\pi}{c_{2k}} \left[\dfrac{LB_k(\lambda - \Delta - \L_0^s, z_{2k+1}(\Nop))}{f_{2k+1}(\Nop)}A - B\right]
	\end{equation*}
	where
	\begin{align*}
		&A = \left(1 - \dfrac{C_{\text{Diag}}}{\pi\sqrt{z_{2k+1}(1)}}\right) \left(1 + \dfrac{1}{f_{2k+1}(1)}\right)^{-1} - \dfrac{C_{\text{Off}}}{\pi K_1(2k+1)^{1+\varepsilon}} \\
		& B = \dfrac{1}{2} \left(1 + \dfrac{1}{f_{2k+1}(1)}\right)^{-1}
	\end{align*}
	which by \eqref{anal_fock_proof_op_est_lower_start_H_2k+1} is also a lower bound for $\H_{2k+1}$. Again, making the constants $K_1$ and $K_2$ in \eqref{anal_fock_def_zk_fk} as large as necessary, we obtain that
	\begin{equation*}
		A \ge 1 - \dfrac{1}{(k+1)^{1+\varepsilon}} \text{  and  } B \le 1 - \dfrac{1}{(k+1)^{1+\varepsilon}} \; \; ,
	\end{equation*}
	which combined with the definition of $c_{2k+1}$ in \eqref{anal_fock_theo_op_estimates_ck} concludes the proof of the lower bound in \eqref{anal_fock_theo_op_estimates_lower_s>1}.
\end{proof}

\begin{proof}[Proof of the upper bound \eqref{anal_fock_theo_op_estimates_upper_s>1}]
	For $k \ge 1$, by the induction hypothesis, we have that
	\begin{equation} \label{anal_fock_proof_op_est_upper_start_H_2k}
		\H_{2k} = \A+^*(\lambda - \Delta - \L_0^s + \H_{2k-1})^{-1} \A+ \le \A+^*(\lambda - \Delta(1 + c_{2k-1}\delta_{2k-1}) - \L_0^s)^{-1} \A+ \; .
	\end{equation}
	As we did before, for every $\psi \in H_n$, we use Lemma \ref{anal_res_CHT_lemma_44} with $S = (\lambda - \Delta(1 + c_{2k-1}\delta_{2k-1}) - \L_0^s)^{-1}$ to separate
	\begin{equation} \label{anal_fock_proof_op_est_upper_diag_+_off}
		\lb\psi, \A+^*(\lambda - \Delta(1 + c_{2k-1}\delta_{2k-1}) - \L_0^s)^{-1} \A+ \psi\rb = \lb\psi, \A+^*S \A+ \psi\rb
	\end{equation}
	into a diagonal and an off-diagonal part, and we treat each separately. For the diagonal part, we apply Lemma \ref{anal_res_CHT_lemma_45}, but this time we want to upper bound
	\begin{align}
		& \int_{\R^2} \dfrac{\hat{V}(q) (\sin \theta)^2 \ud q}{\lambda + |p+q|^2(1 + \frac{c_{2k-1}}{f_{2k}}(LB_{k-1}(\lambda + |p+q|^2 + |p_{1:n}|^{2s} + |q|^{2s},z_{2k}) - f_{2k})) + |p_{1:n}|^{2s} + |q|^{2s}} \nonumber \\
		& \le \frac{f_{2k}}{c_{2k-1}}\int_{\R^2} \dfrac{\hat{V}(q) (\sin \theta)^2 \ud q}{\lambda + |p_{1:n}|^{2s} + |p+q|^2LB_{k-1}(\lambda + |p_{1:n}|^{2s} + |p+q|^2 + |q|^{2s},z_{2k}) + |q|^{2s}} \nonumber \\
		& \le \frac{f_{2k}}{c_{2k-1}}\int_{\R^2} \dfrac{\hat{V}(q) (\sin \theta)^2 \ud q}{\lambda + |p_{1:n}|^{2s} + |p+q|^2LB_{k-1}(\lambda + |p_{1:n}|^{2s} + |p+q|^2,z_{2k})} \; \; , \label{anal_fock_proof_op_est_upper_before_lemma}
	\end{align}
	since $f_{2k-1}(n+1) = f_{2k}(n)$, $z_{2k-1}(n+1) = z_{2k}(n)$. The first inequality is due to the fact that $c_{2k-1} < 1$ and $f_{2k} > 1$ and the second inequality is a consequence of
	\begin{equation*} 
		|p+q|^2LB_{k-1}(\lambda + |p_{1:n}|^{2s} + |p+q|^2,z_{2k}) \le |p+q|^2LB_{k-1}(\lambda + |p_{1:n}|^{2s} + |p+q|^2 + |q|^{2s},z_{2k}) + |q|^{2s}
	\end{equation*}
	Note that the above is equivalent to
	\begin{equation} \label{anal_fock_proof_op_est_upper_trick} 
		LB_{k-1}(\lambda + |p_{1:n}|^{2s} + |p+q|^2,z_{2k}) - LB_{k-1}(\lambda + |p_{1:n}|^{2s} + |p+q|^2 + |q|^{2s},z_{2k}) \le \dfrac{|q|^{2s}}{|p+q|^2} \; .
	\end{equation}
	Thanks to \eqref{appendix_CHT_exp_A4} and the Mean Value Theorem, we have that, for every $k \in \N$ and $x<y \in \R$
	\begin{equation*}
		|LB_{k}(x) - LB_{k}(y)| \le \max_{c \in [x,y]} \dfrac{1}{2(c^2 + c)UB_{k-1}(c,z)}|x-y| \le \frac{1}{x}|x-y|
	\end{equation*}
	The above applied to the difference in \eqref{anal_fock_proof_op_est_upper_trick},  which is positive and hence equals its absolute value, yields
	\begin{equation*}
		LB_{k-1}(\lambda + |p_{1:n}|^{2s} + |p+q|^2,z_{2k}) - LB_{k-1}(\lambda + |p_{1:n}|^{2s} + |p+q|^2 + |q|^{2s},z_{2k}) \le \dfrac{|q|^{2s}}{\lambda + |p_{1:n}|^{2s} + |p+q|^2} \le \dfrac{|q|^{2s}}{|p+q|^2} \; .
	\end{equation*}
	To upper bound the integral in \eqref{anal_fock_proof_op_est_upper_before_lemma} we make use of Lemmas \ref{appendix_CHT_lemma_A2} and \ref{appendix_CHT_lemma_A4} considering $\tilde{\lambda} := \lambda + |p_{1:n}|^{2s}$ instead of $\lambda$, to obtain the upper bound
	\begin{equation*}
		\dfrac{f_{2k}}{c_{2k-1}} \left( \dfrac{\pi}{2}\int_{\lambda + |p_{1:n}|^{2s} + |p|^2}^1 \dfrac{\ud \varrho }{\varrho LB_{k-1}(\varrho,z_{2k})} + \dfrac{C_{\text{Diag}} UB_{k-1}(\lambda + |p_{1:n}|^{2s} + |p|^2)}{\sqrt{z_{2k}}} \right) \; .
	\end{equation*}
	The integral above, by Lemmas \ref{appendix_CHT_lemma_A1} and \ref{appendix_CHT_lemma_A5}, is controlled by
	\begin{align*}
		 \int_{\lambda + |p_{1:n}|^{2s} + |p|^2}^1 \dfrac{\ud \varrho }{\varrho LB_{k-1}(\varrho,z_{2k})}
		& \le \int_{\lambda + |p_{1:n}|^{2s} + |p|^2}^1 \dfrac{\ud \varrho }{(\varrho + \varrho^2) LB_{k-1}(\varrho,z_{2k})} + C\dfrac{UB_{k-1}(\lambda + |p_{1:n}|^{2s} + |p|^2)}{z_{2k}} \\
		& \le 2UB_{k-1}(\lambda + |p_{1:n}|^{2s} + |p|^2) + C\dfrac{UB_{k-1}(\lambda + |p_{1:n}|^{2s} + |p|^2)}{z_{2k}} \; \; .
	\end{align*}
	We deal with the off-diagonal term in the same fashion than in \eqref{anal_fock_proof_op_est_lower_int_2}, upper estimating
	\begin{align*} 
		& n|p| \int_{\R^2} \dfrac{\hat{V}(q)(\sin \theta)^2 \ud q}{[\lambda + |p+q|^2(1 + \frac{c_{2k-1}}{f_{2k}}(LB_{k-1}(\lambda + |p+q|^2 + |p_{1:n}|^{2s} + |q|^{2s},z_{2k}) - f_{2k})) + |p_{1:n}|^{2s} + |q|^{2s}] |p^\prime + q|} \nonumber \\
		& \le n|p| \dfrac{f_{2k}}{c_{2k-1}} \int_{\R^2} \dfrac{\hat{V}(q)(\sin \theta)^2 \ud q}{[\lambda + |p_{1:n}|^{2s} + |p+q|^2LB_{k-1}(\lambda + |p+q|^2 + |p_{1:n}|^{2s} + |q|^{2s},z_{2k})] |p^\prime + q|}
	\end{align*}
	Further, we make use of Lemma \ref{appendix_lemma_2}, this time with $f = LB_{k-1}$ and $g = UB_{k-1}$, to get the upper bound
	\begin{equation*}
		\dfrac{f_{2k}}{c_{2k-1}} \dfrac{n C_{\text{Off}} UB_{k-1}(\lambda + |p_{1:n}|^{2s} + |p|^2, z_{2k})}{z_{2k}}
		\le \dfrac{f_{2k}}{c_{2k-1}} \dfrac{C_{\text{Off}} UB_{k-1}(\lambda + |p_{1:n}|^{2s} + |p|^2, z_{2k})}{K_1(n + 2k)^{1 + 2\varepsilon}} \; \; .
	\end{equation*}
	Putting all the estimates together and noting that $z_{2k}(n) > z_{2k}(1)$, we establish that $\A+^*(\lambda - \Delta(1 + c_{2k-1}\delta_{2k-1}) - \L_0^s)^{-1}\A+$ is upper bounded by
	\begin{equation*}
		\dfrac{\pi}{c_{2k-1}} A^\prime (-\Delta) \delta_{2k} 
	\end{equation*}
	where, by choosing $K_1$ as big as necessary, we obtain
	\begin{equation*}
		A^\prime = 1 + \dfrac{C_{\text{Diag}}}{\pi \sqrt{K_1}(2k)^{1+\varepsilon}} + \dfrac{C}{\pi K_1(2k)^{2+2\varepsilon}} + \dfrac{C_{\text{Off}}}{\pi K_1(2k)^{1+2\varepsilon}} \le 1 + \dfrac{1}{k^{1 + \varepsilon}} \; \; .
	\end{equation*}
	This is enough to see that \eqref{anal_fock_theo_op_estimates_upper_s>1} holds with $c_{2k}$ defined in \eqref{anal_fock_theo_op_estimates_ck}. 
\end{proof}

\section{Proof of \eqref{intro_super_dif_bounds_s>1} in Theorem \ref{intro_theo_main_theo}} \label{section_proof_main_theo_diffusive_bounds}

In this section we finish proving Theorem \ref{intro_theo_main_theo} by using the full power of the iterative estimates provided by Lemma \ref{anal_res_eq_even_odd_inequalities}. This is done by choosing the level of the truncation depending on $\lambda$, i.e., as $\lambda \to 0$, $n \to \infty$ in Lemma \ref{anal_res_eq_even_odd_inequalities}
Again, $C$ denotes a constant, which may change from line to line, but is independent of $p,z,\lambda$ and $k$.

\begin{proof}[Proof of Theorem \ref{intro_theo_main_theo} for $s\ge1$]
	Recall that for $p = (p_1,p_2) \in \R^2$, $\hat{\V}(p) = p_2$ and that $\V \in H_1$ implies that the multiplier of $-\Delta - \L_0^s$ is $|p|^2 + |p|^{2s}$. Let us start with upper bound. By Lemma \ref{anal_res_eq_even_odd_inequalities} and \eqref{set_main_res_laplace_transform_to_resolvent_eq} we get that
	\begin{equation*}
		\dfrac{\lambda^2}{2}D_\V(\lambda) \le \lb \V , \Vs^{2k+1}_1 \rb
		= \lb \V , (\lambda - \Delta - \L_0^s + \H_{2k+1})^{-1} \V\rb \; \; ,
	\end{equation*}
	which by \eqref{anal_fock_theo_op_estimates_lower_s>1} in Theorem \ref{anal_fock_theo_op_estimates} is upper bounded by
		\begin{align}
		& \lb \V, (\lambda - \Delta(1 + c_{2k+1}\delta_{2k+1}) - \L_0^s)^{-1} \V\rb \nonumber\\
		& = \frac{1}{(2\pi)^2} \int_{\R^2} \dfrac{\hat{V}(p)}{|p|^2} \dfrac{|\hat{\V}(p)|^2 \ud p}{\lambda + |p|^2(1 + \frac{c_{2k+1}}{f_{2k+1}} (LB_k(\lambda + |p|^2 + |p|^{2s},z_{2k+1}) - f_{2k+1})) + |p|^{2s}} \nonumber\\
		& \le C\dfrac{f_{2k+1}}{c_{2k+1}}\int_{\R^2} \dfrac{\hat{V}(p) \ud p}{\lambda + |p|^2LB_k(\lambda + |p|^2,z_{2k+1})} \label{proof_of_main_theorem_upper_bound}
	\end{align}
	where we have used \eqref{anal_fock_proof_op_est_upper_trick}. Note that since $\V \in H_1$, the arguments in $f_{2k+1}$ and $z_{2k+1}$ are both $1$ and therefore they are constants which only depend on $k$. 
	Now we conclude exactly as \cite{Cannizaroetal2021}, since the expression above is equal to expression (5.1) in their paper. We include the missing steps in Appendix \ref{section_appendix_end_proof_theo_1} for completeness. 

	Now we proceed to the lower bound. Again, by Lemma \ref{anal_res_eq_even_odd_inequalities} and \eqref{set_main_res_laplace_transform_to_resolvent_eq}, we get that
	\begin{equation*}
		\dfrac{\lambda^2}{2}D_\V(\lambda) \ge \lb \V , \Vs^{2k}_1 \rb
		= \lb \V , (\lambda - \Delta - \L_0^s + \H_{2k})^{-1} \V\rb \; \; ,
	\end{equation*}
	which in turn, by Theorem \ref{anal_fock_theo_op_estimates}, is lower bounded by
	\begin{align}
		\lb \V, (\lambda - \Delta(1 + c_{2k}\delta_{2k}) - \L_0^s)^{-1} \V\rb
		& = \frac{1}{(2\pi)^2} \int_{\R^2} \dfrac{\hat{V}(p)}{|p|^2} \dfrac{|\hat{\V}(p)|^2 \ud p}{\lambda + |p|^2(1 + c_{2k}f_{2k}UB_{k-1}(\lambda + |p|^2 + |p|^{2s},z_{2k})) + |p|^{2s}} \nonumber\\
		& \ge C \int_{\R^2} \dfrac{\hat{V}(p)}{|p|^2} \dfrac{p_2^2 \ud p}{(\lambda + |p|^2)(1 + c_{2k}f_{2k}UB_{k-1}(\lambda + |p|^2,z_{2k})) + |p|^{2s}} \nonumber\\
		&  \ge \dfrac{C}{f_{2k}} \int_{\R^2} \dfrac{\hat{V}(p)}{|p|^2} \dfrac{p_2^2 \ud q}{(\lambda + |p|^2)UB_{k-1}(\lambda + |p|^2,z_{2k}) + |p|^{2s}} \label{proof_of_main_theorem_lower_bound}
	\end{align}
	where 
	we have substituted $c_{2k}$ by its limit as $k \to \infty$ and used the monotonicity of $UB_{k-1}$. 
	Now, note that since all the functions in \eqref{proof_of_main_theorem_lower_bound} but $p \mapsto p_2^2$ are rotationally invariant, the integral has the exact same value as if we replace $p \mapsto p_2^2$ with $p \mapsto p_1^2$. Summing the integrals with $p \mapsto p_2^2$ and $p \mapsto p_1^2$ and diving it by two, we get that expression \eqref{proof_of_main_theorem_lower_bound} is equal to (the $1/2$ is merged into $C$)
		\begin{equation*}
		\dfrac{C}{f_{2k}} \int_{\R^2} \dfrac{\hat{V}(p) \ud p}{(\lambda + |p|^2)UB_{k-1}(\lambda + |p|^2,z_{2k}) + |p|^{2s}} \; \; .
	\end{equation*}
	Thus, an application of \eqref{appendix_lemma_3_final_integral} gives the lower bound
	\begin{align}
		& \dfrac{C}{f_{2k+1}} \left( \int_\lambda^1 \dfrac{\ud \varrho}{\varrho UB_{k-1}(\varrho,z_{2k+1})} - \dfrac{LB_{k}(\lambda, z_{2k+1})}{\sqrt{z_{2k+1}}} \right) \nonumber \\
		& \ge \dfrac{C}{f_{2k+1}} \left( \int_\lambda^1 \dfrac{\ud \varrho}{(\varrho + \varrho^2) UB_{k-1}(\varrho,z_{2k+1})} - \dfrac{LB_{k}(\lambda, z_{2k+1})}{\sqrt{z_{2k+1}}} \right) \nonumber \\
		& \ge \dfrac{C}{f_{2k+1}} \left(LB_{k}(\lambda, z_{2k+1}) - LB_{k}(1, z_{2k+1})  - \dfrac{LB_{k}(\lambda, z_{2k+1})}{\sqrt{z_{2k+1}}} \right) \; \; , \label{proof_main_theorem_lower_final}
	\end{align}
	where the second inequality is a consequence of \eqref{appendix_CHT_exp_A2} in Lemma \ref{appendix_CHT_lemma_A1}. Once again, expression \eqref{proof_main_theorem_lower_final} above reduces to the exact same as the third line in display (5.7) in \cite{Cannizaroetal2021}, and thus we include the end of the proof in Appendix \ref{section_appendix_end_proof_theo_1} for completeness.

	
\end{proof}

\section{Proof of Theorem \ref{intro_theo_invariance_principle}} \label{section_proof_theo_inv_princ}

In this section, we show that our model for $s<1$ is a particular case of the theory developed in Komorowski and Olla \cite{komorowskiandolla2003} of homogenisation for diffusions in divergence free, Gaussian and Markovian random environments. See also Chapters 11 and 12 of the monograph \cite{landimetal2012}. 

Let us consider here the function $\V(\omega) := \omega(0) = (\omega^1(0), \omega^2(0)) = (\V^1(\omega), \V^2(\omega))$. In view of Remark \ref{op_fock_remark_V_in_H_1}, we see that $\V^i \in H_1, i=1,2$. Now, we may write
\begin{equation} \label{proof_theo_inv_princ_X_additive_BM}
	\varepsilon X\left(\frac{t}{\varepsilon^2}\right) = \varepsilon \int_0^{\frac{t}{\varepsilon^2}} \V(\eta_r) \ud r + \sqrt{2}\varepsilon B\left(\frac{t}{\varepsilon^2}\right) \; , \; t\ge0 \; , \varepsilon > 0 \; ,
\end{equation}
and focus on the additive functionals of $(\eta_t)_{t\ge0}$ given by $\int \V^i(\eta_s) \ud s$, for $i=1,2$, since $\varepsilon B(t/\varepsilon^2) \overset{d}{=} B(t)$ for every $\varepsilon>0$ and $t\ge0$. Let $\C := \bigcup_n \oplus_{k \le n} H_k$ be a core for $\L^s$ and $(\L^s)^*$. Let $\mathcal S := (\L^s + (\L^s)^*)/2 =  \L_0^s + \Delta$ be the symmetric part of the generator $\L^s$. For every $\psi \in \C$, let $\|\psi\|_1^2 := \lb \psi, -\L^s \psi \rb = \lb \psi, -\mathcal S \psi \rb$ be a norm and 
\begin{equation} \label{proof_theo_inv_princ_def_H_{-1}_norm}
	\|\psi\|_{-1}^2 := \lim\limits_{\lambda \to 0} \lb \psi, (\lambda - \mathcal S)^{-1} \psi \rb
\end{equation}
be another norm.  By \cite[Theorem 2.2]{sethuramanetal2000}, for every $t\ge0$, it holds that
\begin{equation} \label{proof_theo_inv_princ_H_-1_norm_bound_variance}
	\E\left[\sup_{0\le t^\prime \le t} \left|\int_0^{t^\prime} \mathcal V^i(\eta_r) \ud r\right|^2\right] \le Ct\|\V^i\|_{-1}  \le Ct \; , i=1,2 \; ,
\end{equation}
where the last inequality is a consequence of
\begin{equation}
	\|\V^i\|_{-1}^2 = \frac{1}{(2\pi)^2} \int_{\R^2} \dfrac{\hat{V}(p)}{|p|^2} \dfrac{|\hat{\V^i}(p)|^2 \ud p}{|p|^2 + |p|^{2s}}
	\le C \int_{\R^2} \dfrac{\hat{V}(p) \ud p}{|p|^{2s}}
	\le C \int_0^1 r^{1-2s} \ud r \le C \text{ , for } i=1,2
\end{equation}
since $s<1$, as discussed previously in \eqref{anal_res_remark_level_1_eq} in Remark \ref{anal_res_remark_level_1} for $i=1$.
Note that \eqref{proof_theo_inv_princ_H_-1_norm_bound_variance} proves the upper bound \eqref{intro_dif_bounds_s<1} in Theorem \ref{intro_theo_invariance_principle}. The lower bound follows from the Yaglom-reversibility \eqref{set_main_res_yaglom_reversibility}.
So now we show that our model, for $s<1$, is a particular case of the general framework of divergence-free, Gaussian and Markovian environments treated in \cite[Section 6]{komorowskiandolla2003}.

\begin{proof}[Proof of Theorem \ref{intro_theo_invariance_principle}]
In Section 6 of \cite{komorowskiandolla2003}, the same SDE as in \eqref{intro_def_main_SDE} is considered, with a dynamic random environment $(\omega_t)_{t\ge0}$ which is divergence-free, Gaussian and Markovian. Moreover, they assume that, in $d=2$, the space-time correlations of the drift field $\omega$ satisfy 
expression (1.2) in page $181$, which reads as
\begin{equation} \label{proof_theo_inv_princ_covariance_of_field}
	R(t,x) = \int_{\R^2} e^{\iota x\cdot p} \exp\{-|p|^{2\beta}t\} \frac{a(|p|)}{|p|^{2\alpha}} \left(\mathbf I - \frac{p \otimes p}{|p|^2}\right) \ud p = \int_{\R^2} e^{\iota x\cdot p} \exp\{-|p|^{2\beta}t\} \frac{a(|p|)}{|p|^{2\alpha+2}} \left(\mathbf I |p|^2- p \otimes p\right) \ud p \; ,
\end{equation}
where $a : [0,\infty) \to [0,\infty)$ is a compactly supported and bounded cut-off function, $\beta \ge 0$ and $\alpha < 1$. Also, the notation $p \otimes p$ represents the canonical tensor product in $\R^2$ and $\mathbf I$ the identity $2 \times 2$ matrix. Since here we consider the dyamics in \eqref{intro_def_environmental_process_omega}, we identify $\beta$ in \eqref{proof_theo_inv_princ_covariance_of_field} with $s$. Also, since $\hat{V}$ is rotationally invariant and has compact support, we may identify $a(|p|)$ in \eqref{proof_theo_inv_princ_covariance_of_field} with $\hat{V}(p)$. Now, note that, for $p = (p_1,p_2) \in \R^2$,
\begin{equation*}
	\left(\mathbf I |p|^2- p \otimes p\right) = \begin{pmatrix}
		p_1^2 + p_2^2 & 0\\
		0 & p_1^2 + p_2^2
	\end{pmatrix}
	-
	\begin{pmatrix}
		p_1p_1 & p_1p_2\\
		p_2p_1 & p_2p_2
	\end{pmatrix}
	=
	\begin{pmatrix}
		p_2^2 & -p_1p_2\\
		-p_2p_1 & p_1^2 \; \; ,
	\end{pmatrix} \; \; .
\end{equation*}
In view of \eqref{op_fock_def_brackets} and \eqref{op_fock_convention_fourier}, we get that, for every $\psi^j \in H_1$, $j=1,2$, given by $\psi^j(\omega) = \int_{\R^2} f^j_1(x) \omega^1(x) \ud x + \int_{\R^2} f^j_2(x) \omega^2(x) \ud x$, (in what follows we suppress $p$ from $\hat{f^i_j}(p)$)
\begin{align*}
	\lb \psi^1, \psi^2 \rb & = \frac{1}{(2\pi)^2} \int_{\R^2} \frac{\widehat{V}(p)}{|p|^2} \widehat{\psi^1}(p) \overline{\widehat{\psi^2}(p)} \ud p \\
	& = \frac{1}{(2\pi)^2} \int_{\R^2} \frac{\widehat{V}(p)}{|p|^2} [p_2^2\hat{f_1^1}\overline{\hat{f_1^2}} - p_2p_1\hat{f_1^1}\overline{\hat{f_2^2}} - p_1p_2\hat{f_2^1}\overline{\hat{f_1^2}} + p_1^2\hat{f_2^1}\overline{\hat{f_2^2}}] \ud p \\
	& = \frac{1}{(2\pi)^2} \int_{\R^2} \frac{\widehat{V}(p)}{|p|^2}
	\begin{pmatrix}
		\hat{f_1^1} & \hat{f_2^1}
	\end{pmatrix}
	\begin{pmatrix}
		p_2^2 & -p_1p_2\\
		-p_2p_1 & p_1^2 \; \; ,
	\end{pmatrix}
	\overline{\begin{pmatrix}
		\hat{f_1^2}\\
		\hat{f_2^2}
	\end{pmatrix}}
	\ud p \; \; .
\end{align*}
With this observation, we see that $\alpha < 1$ in \eqref{proof_theo_inv_princ_covariance_of_field} translates to $\alpha = 0$. With the same argument, we conclude that the law of $\cGFF$ satisfies assumption (E) in Section 6 of \cite{komorowskiandolla2003} with $\alpha = 0$. Therefore, since $s = \beta < 1$, by
\cite[Theorem 6.3]{komorowskiandolla2003}, we get Theorem \ref{intro_theo_invariance_principle} with the covariance matrix $D$ given by
\begin{equation} \label{proof_theo_inv_principle_covariance_matrix_D}
	D = 2[\delta_{i,j} + \lb \psi_*^i, \psi_*^j \rb_1]_{i,j} \; , \; i,j=1,2 \; ,
\end{equation}
where the objects $\psi_*^i$, for $i=1,2$ satisfy $\lim_{\lambda \downarrow 0} \|\psi_\lambda^i - \psi_*^i\|_1 = 0$ for $\psi_\lambda^i$ solution to the resolvent equations
\begin{equation*}
	\lambda \psi_\lambda^i - \L^s\psi_\lambda^i = - \V^i \; , \;  \lambda > 0 \; .
\end{equation*}
The inner product $\lb \cdot , \cdot , \rb_1$ is defined through polarisation by
\begin{equation*}
	\lb \psi_*^i, \psi_*^j \rb_1 := \frac{1}{4} \left(\|\psi_*^i + \psi_*^j\|_1^2 + \|\psi_*^i - \psi_*^j\|_1^2 \right) \; , \; i,j=1,2 \; \; .
\end{equation*}

\end{proof}

\section{Proof of Theorem \ref{intro_theo_log_div}} \label{section_proof_theo_log_div}

In this section we prove Theorem \ref{intro_theo_log_div} by making use of the first upper and lower bounds provided by Lemma \ref{anal_res_eq_even_odd_inequalities}, i.e., the estimates obtained for $n=1$.
When $\omega_t = \omega_t^\gamma$ is the solution to \eqref{intro_def_environmental_process_omega_log_div}, the dominant terms in the estimates are once again the ones coming from the dynamics of the environment, as in the case of $s<1$ in Theorem \ref{intro_theo_main_theo}. This is the reason why we can find matching upper and lower bounds just going to the first two estimates.

Since $(-\Delta)$ is a self-adjoint, positive operador, we can make sense of the operator $(\log(e + (-\Delta)^{-1}))^\gamma$ for every $\gamma > 0$ through 
its Fourier multiplier, in the spirit of Proposition \ref{op_fock_prop_invariant_measure_fSHE}, given by
\begin{equation} \label{proof_log_div_fourier_multiplier}
	\sigma^\gamma(p_{1:n}) = \sum_{i=0}^{n} |p_i|^2(\log(e + |p_i|^{-2}))^\gamma \; \; .
\end{equation}
Expression \eqref{proof_log_div_fourier_multiplier} is associated with the generator $\L_0^\gamma$ of the process $(\omega_t^\gamma)_{t\ge0}$ solution to \eqref{intro_def_environmental_process_omega_log_div}. Therefore, since we have the correction $(\log(e + (-\Delta)^{-1}))^\frac{\gamma}{2}$ in front of the noise in \eqref{intro_def_environmental_process_omega_log_div}, Proposition \ref{op_fock_prop_invariant_measure_fSHE} holds true with $-(-\Delta)^s$ replaced by $(\log(e + (-\Delta)^{-1}))^\gamma\Delta$ and thus the dynamics in \eqref{intro_def_environmental_process_omega_log_div} preserves the law of $\cGFF$ as invariant measure for every $\gamma > 0$. Note also that $(\log(e + x^{-1}))^\gamma > 1$ for every $x \ge 0$.

Let us start with some calculus, which are analogous results to the ones in Lemma \ref{appendix_CHT_lemma_A1}. For every $\gamma>0$, $\gamma \neq 1$, the following holds:
\begin{align}
	& \partial_x (\log (e + x^{-1}))^\gamma = - \frac{\gamma}{2} \frac{1}{(ex^2 + x)(\log (e + x^{-1}))^{1-\gamma}} 
	\ge -\gamma \frac{1}{x} \; \; , \label{proof_theo_log_div_deriv_gamma}\\
	& \partial_x (\log (e + x^{-1}))^{1-\gamma} = - \frac{1-\gamma}{2} \frac{1}{(ex^2 + x)(\log (e + x^{-1}))^{\gamma}} \; \; .\label{proof_theo_log_div_deriv_1-gamma}
\end{align}

\begin{proof}[Proof of Theorem \ref{intro_theo_log_div}]
	We denote $\L^\gamma = \L_0^\gamma + \A+ - \A+^* + \Delta$. Here again we only consider $\V := \V^1 \in H_1$ given by $\V(\omega) := \omega^1(0)$, as in Theorem \ref{intro_theo_main_theo}. Recall that $\hat{\V}(p) = p_2$ for $p = (p_1,p_2)$. Thus, by taking $n=1$ in Lemma \ref{anal_res_eq_even_odd_inequalities}, we arrive at
	\begin{align} 
		\lb \V , (\lambda - \L^\gamma)^{-1}\V \rb \le \lb \V , (\lambda - \Delta - \L_0^\gamma)^{-1}\V \rb
		=& \frac{1}{(2\pi)^2} \int_{\R^2} \dfrac{\hat{V}(p)}{|p|^2} \dfrac{|\hat{\V}(p)|^2 \ud p}{\lambda + |p|^2 + |p|^2(\log(e + |p|^{-2}))^\gamma} \nonumber \\
		& \le C \int_{\R^2} \dfrac{\hat{V}(p) \ud p}{\lambda + |p|^2(\log(e + (\lambda + |p|^2)^{-1}))^\gamma} \label{proof_theo_log_div_upper_bound}
	\end{align}
	Now, adapting \eqref{appendix_CHT_exp_A9} in Lemma \ref{appendix_CHT_lemma_A2} and \eqref{appendix_lemma_3_squeeze_r^2_in} in Lemma \ref{appendix_lemma_3}, we have that expression \eqref{proof_theo_log_div_upper_bound} is upper bounded by
	\begin{equation} \label{proof_theo_log_div_integral_computed}
		C\int_\lambda^1 \frac{\ud s}{s(\log(e + s^{-1}))^\gamma} + C \le C\int_\lambda^1 \frac{\ud s}{(es^2 + s)(\log(e + s^{-1}))^\gamma} + C
	\end{equation}
	Therefore, by \eqref{proof_theo_log_div_deriv_1-gamma} and \eqref{set_main_res_laplace_transform_to_resolvent_eq}, we see that
	\begin{equation}
		D_\V(\lambda) \le C (\log (e + \lambda^{-1}))^{1-\gamma} \overset{\lambda \to 0}{\le} C|\log \lambda|^{1-\gamma} \; \; . \label{proof_theo_log_div_upper_bound_first_level}
	\end{equation}
Note that if $\gamma > 1$, \eqref{proof_theo_log_div_upper_bound_first_level} is bounded by a constant and this is enough to show the diffusive bounds in \eqref{intro_dif_bounds_gamma>1} by the Yaglom-reversibility \eqref{set_main_res_yaglom_reversibility}. Also, if $\gamma = 1$ in \eqref{proof_theo_log_div_integral_computed}, then by \eqref{appendix_CHT_exp_A2} with $k=0$, we get the upper bound in \eqref{intro_loglog_bounds_gamma_1}.

Now, let us proceed to the lower bound by also taking $n=1$ in Lemma \ref{anal_res_eq_even_odd_inequalities}, for the case $\gamma \in [\frac{1}{2},1]$. First, note that by adapting Lemma \ref{appendix_lemma_2}, the off diagonal term in the first lower estimate is bounded above by a constant $C$. Also, note that by adapting \eqref{appendix_CHT_exp_A8} in Lemma \ref{appendix_CHT_lemma_A2} and using \eqref{proof_theo_log_div_deriv_1-gamma} again, we see that, for a constant $D > 0$
\begin{align} 
	&\lb \V , (\lambda - \L^\gamma)^{-1}\V \rb \ge
	\lb \V , (\lambda - \Delta - \L_0^\gamma + \A+^*(\lambda - \Delta - \L_0^\gamma)\A+)^{-1}\V \rb \nonumber \\
	& \ge C \int_{\R^2} \dfrac{\hat{V}(p)}{|p|^2} \dfrac{|\hat{\V}(p)|^2 \ud p}{\lambda + |p|^2(\log(e + |p|^{-2}))^\gamma + D|p|^2(1 + (\log (e + (\lambda + |p|^2)^{-1}))^{1-\gamma})} \nonumber \\
	& \ge C \int_{\R^2}  \dfrac{\hat{V}(p) \ud p}{\lambda + |p|^2(\log(e + |p|^{-2}))^\gamma + D|p|^2(1 + (\log (e + (\lambda + |p|^2)^{-1}))^{1-\gamma})}  \nonumber \\
	& \ge C \int_{\R^2}  \dfrac{\hat{V}(p) \ud p}{\lambda + |p|^2(\log(e + |p|^{-2}))^\gamma}
	\ge C \int_{\R^2}  \dfrac{\hat{V}(p) \ud p}{\lambda + |p|^2(\log(e + (\lambda + |p|^2)^{-1}))^\gamma} \label{proof_theo_log_div_lower_bound} \; \; .
\end{align}
The third inequality is a result of the same argument as in \eqref{proof_of_main_theorem_lower_bound}, the fourth inequality is true because $\gamma \in [\frac{1}{2},1] \Rightarrow 1-\gamma \le \gamma$ and thus we may absorb the lower order terms into $|p|^2(\log(e + |p|^{-2}))^\gamma$ by changing the constant $C$. The fifth inequality is due to
an application of the Mean Value Theorem together with the inequality in \eqref{proof_theo_log_div_deriv_gamma}, in the same spirit of \eqref{anal_fock_proof_op_est_upper_trick}. Once again, by adapting \eqref{appendix_CHT_exp_A9} in Lemma \ref{appendix_CHT_lemma_A2}, we get that \eqref{proof_theo_log_div_lower_bound} is lower bounded by
\begin{equation*}
	C\int_\lambda^1 \frac{\ud s}{s(\log(e + s^{-1}))^\gamma} - C \ge C\int_\lambda^1 \frac{\ud s}{(es^2 + s)(\log(e + s^{-1}))^\gamma} - C \; \; .
\end{equation*}
Therefore, if $\gamma \in [\frac{1}{2},1)$, by \eqref{proof_theo_log_div_deriv_1-gamma} and \eqref{set_main_res_laplace_transform_to_resolvent_eq}, we get that
\begin{equation} \label{proof_theo_log_div_lower_bound_first_level}
	D_\V(\lambda) \ge C (\log (e + \lambda^{-1}))^{1-\gamma} - C \ge C (\log (e + \lambda^{-1}))^{1-\gamma} \overset{\lambda \to 0}{\ge} C|\log \lambda|^{1-\gamma} \; \; ,
\end{equation}
and if $\gamma = 1$, by \eqref{appendix_CHT_exp_A2} with $k=0$, we get the lower bound in \eqref{intro_loglog_bounds_gamma_1},
which concludes the proof of Theorem \ref{intro_theo_log_div}.
\end{proof}

\appendix

\section{Technical lemmas I} \label{section_appendix_their_lemmas}

In this section, for completeness, we list some important technical lemmas used throughout the estimates in the proofs of Theorem \ref{anal_fock_theo_op_estimates} and Theorem \ref{intro_theo_main_theo}, all of them due to Toninelli et al. \cite{Cannizaroetal2021}.

\begin{lemma} \label{appendix_CHT_lemma_A1}
	For $k \in \N$ let $L$, $LB_k$ and $UB_k$ be the functions defined in \eqref{anal_fock_def_L} and \eqref{anal_fock_def_LB_UB}. Then, the three are decreasing in the first variable and increasing in the second. For every $x > 0$ and $z \ge 1$, the following holds true
	\begin{equation} \label{appendix_CHT_exp_A1}
		\begin{split}
			& 1 \le LB_k(x,z) \le \sqrt{L(x,z)} \; \; , \\
			& 1 \le \sqrt{z} \le \sqrt{L(x,z)} \le UB_k(x,z) \le L(x,z) \; \; .
		\end{split}
	\end{equation}
	Furthermore, for every $0 < a < b$, one has
	\begin{align}
		& \int_a^b \dfrac{\ud x}{(x^2 + x)UB_k(x,z)} = 2(LB_{k+1}(a,z) - LB_{k+1}(b,z)) \; \; , \label{appendix_CHT_exp_A2} \\
		& \int_a^b \dfrac{\ud x}{(x^2 + x)LB_k(x,z)} \le 2(UB_{k}(a,z) - UB_{k}(b,z)) \; \; . \label{appendix_CHT_exp_A3}
	\end{align}
	Finally, it also holds that
	\begin{equation} \label{appendix_CHT_exp_A4}
		\begin{split}
			& \partial_xL(x,z) = -\dfrac{1}{x^2 + x} \; \; , \quad
			 \partial_xLB_k(x,z) = -\dfrac{1}{2(x^2 + x)UB_{k-1}(x,z)} \; \; , \\
			 & \partial_xUB_k(x,z) = -\dfrac{1}{2(x^2 + x)LB_{k}(x,z)}\left(1 + \dfrac{(\frac{1}{2} \log L(x,z))^k}{k! LB_k(x,z)}\right) \; \; .
		\end{split}
	\end{equation}
\end{lemma}

\begin{lemma} \label{appendix_CHT_lemma_A2}
	Let $V$ be as in \eqref{intro_def_mollifier_V}. Let $z > 1$ and $f(\cdot, z) : [0,\infty) \mapsto [1,\infty)$ be a strictly decreasing and differentiable function, such that 
	\begin{equation} \label{appendix_CHT_exp_A7}
		-\dfrac{f(x)}{x} \le f^\prime(x) < 0 \text{ for all } x \in \R
	\end{equation} 
	and the function $g(\cdot, z) : [0,\infty) \mapsto [1,\infty)$ a strictly decreasing function such that $g(x,z)f(x,z) \ge z$. Then, there exists a constant $C_{\text{Diag}} > 0$ such that, for all $z > 1$, one gets the bound
	\begin{equation} \label{appendix_CHT_exp_A8}
		\left|\int_{\R^2}\dfrac{\hat{V}(q)(\sin \theta)^2 \ud q}{\lambda + |p + q|^2f(\lambda + |p + q|^2,z)} - \dfrac{\pi}{2}\int_{\lambda + |p|^2}^1 \dfrac{\ud \varrho}{\varrho f(\varrho,z)}\right| \le C_{\text{Diag}} \dfrac{g(\lambda + |p|^2,z)}{\sqrt{z}}
	\end{equation}
	where $p = \sum_{i=1}^{n} p_i$ for some $n \in \N$ and $p_1,\dots, p_n \in \R^2$ and $\theta$ is the angle between $p$ and $q$. The second integral is zero if $\lambda + |p|^2 \ge 1$.
	Moreover, for $\lambda \le 1$,
	\begin{equation} \label{appendix_CHT_exp_A9}
		\left|\dfrac{1}{2}\int_{\R^2}\dfrac{\hat{V}(q) \ud q}{\lambda + |q|^2f(\lambda + |q|^2,z)} - \dfrac{\pi}{2}\int_{\lambda}^1 \dfrac{\ud \varrho}{\varrho f(\varrho,z)}\right| \le C_{\text{Diag}} \dfrac{g(\lambda,z)}{\sqrt{z}}
	\end{equation}
\end{lemma}

\begin{lemma} \label{appendix_CHT_lemma_A4}
	The functions $UB_k(\cdot, z)$ and $LB_k(\cdot, z)$ satisfy the conditions of the previous lemmas. 
\end{lemma}

\begin{lemma} \label{appendix_CHT_lemma_A5}
	For every $z\ge1$, $\lambda \in \R_+$ and $p \in \R^2$ such that $\lambda + |p|^2 \le 1$, one has
	\begin{equation*}
		\left|\int_{\lambda + |p|^2}^1 \dfrac{\ud \varrho}{\varrho LB_k(\varrho,z)}
		 - \int_{\lambda + |p|^2}^1 \dfrac{\ud \varrho}{(\varrho + \varrho^2) LB_k(\varrho,z)}\right| \le \dfrac{UB_k(\lambda + |p|^2,z)}{z} \; \; .
	\end{equation*}
\end{lemma}

\section{Technical lemmas II} \label{section_appendix_our_lemmas}

The two lemmas in this section are modifications of Lemmas \ref{appendix_CHT_lemma_A2} and Lemma A.3 in Toninelli et al. \cite{Cannizaroetal2021}. Throughout this section we use a generic constant $C$ which may change from line to line, but is always independent of $p, q, z, \lambda, k$ and $n$.

\begin{lemma} \label{appendix_lemma_3}
	Let $s\ge1$ in \eqref{intro_def_environmental_process_omega} and $\tilde{\lambda} := \lambda + |p_{1:n}|^{2s}$. Then, there exists a constant $C_{\text{Diag}} > 0$ such that, for all $z > 1$ and every $k \ge 0$, we get the bound
	\begin{equation} \label{appendix_lemma_3_main_recursive}
		\left|\int_{\R^2}\dfrac{\hat{V}(q)(\sin \theta)^2 \ud q}{(\tilde{\lambda} + |p + q|^2)UB_k(\tilde{\lambda} + |p + q|^2,z) + |q|^{2s}} - \dfrac{\pi}{2}\int_{\tilde{\lambda} + |p|^2}^1 \dfrac{\ud \varrho}{\varrho UB_k(\varrho,z)}\right| \le C_{\text{Diag}} \dfrac{LB_{k+1}(\tilde{\lambda} + |p|^2,z)}{\sqrt{z}}
	\end{equation}
	where $p = \sum_{i=1}^{n} p_i$ for some $n \in \N$ and $p_1,\dots, p_n \in \R^2$ and $\theta$ is the angle between $p$ and $q$. The second integral is zero if $\tilde{\lambda} + |p|^2 \ge 1$.
	Moreover, for $\lambda \le 1$,
	\begin{equation} \label{appendix_lemma_3_final_integral}
		\left|\dfrac{1}{2}\int_{\R^2}\dfrac{\hat{V}(q) \ud q}{(\lambda + |q|^2)UB_k(\lambda + |q|^2,z) + |q|^{2s}} - \dfrac{\pi}{2}\int_{\lambda}^1 \dfrac{\ud \varrho}{\varrho UB_k(\varrho,z)}\right| \le C_{\text{Diag}} \dfrac{LB_{k+1}(\lambda,z)}{\sqrt{z}}
	\end{equation}
\end{lemma}

\begin{proof}
	Since $z$ is fixed, we suppress the dependence of $UB_k$ and $LB_k$ on it. A fact used multiple times here is that for all $a,b > 0$ and $z \ge 1$ we have
	\begin{equation} \label{appendix_ineq_lemma_UB_and_LB}
		\dfrac{1}{UB_k(a + b,z)} \le \dfrac{LB_k(a+b,z)}{z} \le \dfrac{LB_k(a,z)}{z} \le \dfrac{LB_k(a,z)}{\sqrt{z}} \le \dfrac{LB_{k+1}(a,z)}{\sqrt{z}}\; .
	\end{equation}
	First, we separate the left hand side of \eqref{appendix_lemma_3_main_recursive} into three terms
	\begin{align}
		& \left|\int_{\R^2}\dfrac{\hat{V}(q)(\sin \theta)^2 \ud q}{(\tilde{\lambda} + |p + q|^2)UB_k(\tilde{\lambda} + |p + q|^2) + |q|^{2s}} -
		\int_{\R^2}\dfrac{\hat{V}(q)(\sin \theta)^2 \ud q}{(\tilde{\lambda} + |p|^2 + |q|^2)UB_k(\tilde{\lambda} + |p|^2 + |q|^2) + |q|^{2s}}\right| \label{appendix_lemma_3_first_term} \\
		& + \left|\int_{\R^2}\dfrac{\hat{V}(q)(\sin \theta)^2 \ud q}{(\tilde{\lambda} + |p|^2 + |q|^2)UB_k(\tilde{\lambda} + |p|^2 + |q|^2) + |q|^{2s}} -
		\int_{\R^2}\dfrac{\hat{V}(q)(\sin \theta)^2 \ud q}{(\tilde{\lambda} + |p|^2 + |q|^2)UB_k(\tilde{\lambda} + |p|^2 + |q|^2)}\right| \label{appendix_lemma_3_second_term} \\
		& + \left|\int_{\R^2}\dfrac{\hat{V}(q)(\sin \theta)^2 \ud q}{(\tilde{\lambda} + |p|^2 + |q|^2)UB_k(\tilde{\lambda} + |p|^2 + |q|^2)} -
		\dfrac{\pi}{2}\int_{\tilde{\lambda} + |p|^2}^1 \dfrac{\ud \varrho}{\varrho UB_k(\varrho)}\right| \label{appendix_lemma_3_third_term}
	\end{align}
	Note that \eqref{appendix_lemma_3_first_term} and \eqref{appendix_lemma_3_third_term} have the same flavour as (A.11) and (A.12) in \cite[Lemma A.2]{Cannizaroetal2021}, 
	respectively.
	In fact, we handle them almost indentically 
	and we add the proof here for completeness. The main difference is then in the term \eqref{appendix_lemma_3_second_term}. We start with \eqref{appendix_lemma_3_first_term}. Note that under the restriction $|p+q| < |p|$, we may bound each integral individually. In fact, for the first we use $(\sin \theta)^2 \le \frac{|p+q|^2}{|p|^2}$ to get
	\begin{align*}
		\int_{|p+q| < |p|}\dfrac{\hat{V}(q)(\sin \theta)^2 \ud q}{(\tilde{\lambda} + |p + q|^2)UB_k(\tilde{\lambda} + |p + q|^2) + |q|^{2s}} & \le
		|p|^{-2} \int_{|p+q| < |p|}\dfrac{\hat{V}(q) \ud q}{UB_k(\tilde{\lambda} + |p + q|^2)}\\
		& \le \frac{|p|^{-2}}{UB_k(\tilde{\lambda} + |p|^2)} \int_{|p+q| < |p|} \ud q
		\le C\frac{LB_{k+1}(\tilde{\lambda} + |p|^2)}{\sqrt{z}} \; \; .
	\end{align*}
	For the second, we see that $|p+q| < |p| \Rightarrow |q| < 2|p|$ and therefore $(\tilde{\lambda} + |p|^2 + |q|^2)UB_k(\tilde{\lambda} + |p|^2 + |q|^2) + |q|^{2s} \ge |p|^2UB_k(\tilde{\lambda} + 5|p|^2)$ and again $\int_{|p+q| < |p|} \ud q \le C|p|^2$.
	
	For the region $|p+q| \ge |p|$, let $h(x) = xUB_k(x)$, which by Lemma \ref{appendix_CHT_lemma_A4} satisfies \eqref{appendix_CHT_exp_A7} and thus $|h^\prime(x)| \le 2|f(x)|$. So by the Mean Value Theorem,
	\begin{equation} \label{appendix_lemma_1_aux_function_h_new}
		|h(x) - h(y)| \le 2|x-y| UB_k(\min\{x,y\}) \text{  and  } |p+q|^2 - |p|^2 - |q|^2 = |p||q| \cos \theta \; .
	\end{equation} Note that by trashing both $|q|^{2s} \ge 0$ in the denominator of the difference in \eqref{appendix_lemma_3_first_term}, over $|p+q| \ge |p|$ the difference is bounded by
	\begin{align*}
		& \int_{|p+q| \ge |p|} \dfrac{\hat{V}(q)(\sin \theta)^2 |h(\tilde{\lambda} +  |p+q|^2) - h(\tilde{\lambda} + |p|^2 + |q|^2)| \ud q}{(\tilde{\lambda} + |p+q|^2)(\tilde{\lambda} + |p|^2 + |q|^2)UB_k(\tilde{\lambda} + |p+q|^2)UB_k(\tilde{\lambda} + |p|^2 + |q|^2)} \\
		& \le C\int_{|p+q| \ge |p|} \dfrac{\hat{V}(q)(\sin \theta)^2 |p||q||\cos \theta| \ud q}{(\tilde{\lambda} + |p+q|^2)(\tilde{\lambda} + |p|^2 + |q|^2)UB_k(\tilde{\lambda} + \max\{|p+q|^2, |p|^2 + |q|^2\})} \\
		& \le C|p|\int_{|p+q| \ge |p|} \dfrac{\hat{V}(q)|q| \ud q}{(\tilde{\lambda} + |p+q|^2)(\tilde{\lambda} + |p|^2 + |q|^2)UB_k(\tilde{\lambda} + 2|p|^2 + 2|q|^2)} \\
		& \le C\frac{LB_k(\tilde{\lambda} + |p|^2)}{\sqrt{z}}|p|\int_{|p+q| \ge |p|} \dfrac{\hat{V}(q)|q| \ud q}{(\tilde{\lambda} + |p+q|^2)(\tilde{\lambda} + |p|^2 + |q|^2)} \le
		C\frac{LB_{k+1}(\tilde{\lambda} + |p|^2)}{\sqrt{z}} \; \; ,
	\end{align*}
	where the last inequality is a consequence of the integral in the last line being of order $|p|^{-1}$. To see that, further divide the integral into the regions $|q| \ge 2|p|$ and $|q| < 2|p|$. For the first, note that $|q| \ge 2|p| \Rightarrow |p+q| \ge \frac{|q|}{2}$
	\begin{equation*}
		\int_{|p+q| \ge |p| \atop |q| \ge 2|p|} \dfrac{\hat{V}(q)|q| \ud q}{(\tilde{\lambda} + |p+q|^2)(\tilde{\lambda} + |p|^2 + |q|^2)}
		\le C \int_{|q| \ge |p|} \hat{V}(q)|q|^{-3} \ud q \le C \int_{|p|}^1 r^{-2} \ud r \le \dfrac{C}{|p|} \; ,
	\end{equation*}
	while for the second
	\begin{equation*}
		\int_{|p+q| \ge |p| \atop |q| < 2|p|} \dfrac{\hat{V}(q)|q| \ud q}{(\tilde{\lambda} + |p+q|^2)(\tilde{\lambda} + |p|^2 + |q|^2)}
		\le \dfrac{C}{(\tilde{\lambda} + |p|^2)^2} \int_{|q| < 2|p|} |q| \ud q \le \dfrac{C}{|p|^4} \int_0^{2|p|} r^2 \ud r \le \dfrac{C}{|p|} \; .
	\end{equation*}
	
	This concludes the estimate of the first term.
	
	Now, we move to \eqref{appendix_lemma_3_third_term} and conclude with \eqref{appendix_lemma_3_second_term} at the end since we will need \eqref{appendix_lemma_3_third_term} for its proof. Here, we consider the first integral over the region $|q|^2 \ge 1 - (\tilde{\lambda} + |p|^2)$, which implies $(\tilde{\lambda} + |p|^2 + |q|^2)^{-1} \le 1$. Using \eqref{appendix_ineq_lemma_UB_and_LB}, we obtain the following upper bound
	\begin{equation*}
		\int_{\R^2}\dfrac{\hat{V}(q)(\sin \theta)^2 \ud q}{(\tilde{\lambda} + |p|^2 + |q|^2)UB_k(\tilde{\lambda} + |p|^2 + |q|^2)} \le 
		\dfrac{LB_k(\tilde{\lambda} + |p|^2)}{z} \int_{\R^2} \hat{V}(q) \ud q \le C\dfrac{LB_{k+1}(\tilde{\lambda} + |p|^2)}{\sqrt{z}} \; .
	\end{equation*}
	Still the first integral in \eqref{appendix_lemma_3_third_term} but now in the complement of the previous region, we first observe that since $\hat{V}$ is smooth and rotationally invariant, there exists a constant $C > 0$ such that $|\hat{V}(q) - \hat{V}(0)| < C|q|^2$ for $|q| \le 1$. Then, we may re-write that integral as
	\begin{align}
		& \int_{|q|^2 < 1 - (\tilde{\lambda} + |p|^2)} \dfrac{\hat{V}(0)(\sin \theta)^2 \ud q}{(\tilde{\lambda} + |p|^2 + |q|^2)UB_k(\tilde{\lambda} + |p|^2 + |q|^2)} \label{appendix_proof_lemma_1_int_4_new} \\
		& + \int_{|q|^2 < 1 - (\tilde{\lambda} + |p|^2)} \dfrac{(\hat{V}(q) - \hat{V}(0))(\sin \theta)^2 \ud q}{(\tilde{\lambda} + |p|^2 + |q|^2)UB_k(\tilde{\lambda} + |p|^2 + |q|^2)} \; . \label{appendix_proof_lemma_1_int_5_new}
	\end{align}
	Passing the integral in \eqref{appendix_proof_lemma_1_int_4_new} into polar coordinates and then setting $s = \tilde{\lambda} + |p|^2 + r^2$, we get
	\begin{equation} \label{appendix_proof_lemma_1_change_of_variables}
		\int_0^{2\pi} (\sin \theta)^2\ud \theta \int_0^{\sqrt{1 - \tilde{\lambda} - |p|^2}} \dfrac{ r \ud r}{(\tilde{\lambda} + |p|^2 + r^2)UB_k(\tilde{\lambda} + |p|^2 + r^2)} = \dfrac{\pi}{2} \int_{\tilde{\lambda} + |p|^2}^1 \dfrac{\ud s}{sUB_k(s)}
	\end{equation}
	Lastly, we control the integral in \eqref{appendix_proof_lemma_1_int_5_new} using $|\hat{V}(q) - \hat{V}(0)| < C|q|^2$ for $|q| \le 1$ and \eqref{appendix_ineq_lemma_UB_and_LB}
	\begin{align*}
		& \int_{|q|^2 < 1 - (\tilde{\lambda} + |p|^2)} \dfrac{|\hat{V}(q) - \hat{V}(0)|(\sin \theta)^2 \ud q}{(\tilde{\lambda} + |p|^2 + |q|^2)UB_k(\tilde{\lambda} + |p|^2 + |q|^2)} \\
		& \le C\dfrac{LB_k(\tilde{\lambda} + |p|^2)}{z}\int_{|q|^2 < 1 - (\tilde{\lambda} + |p|^2)} \dfrac{|q|^2 \ud q}{\tilde{\lambda} + |p|^2 + |q|^2}
		\le C\dfrac{LB_{k+1}(\tilde{\lambda} + |p|^2)}{\sqrt{z}}\int_{|q| < 1} \ud q
		\le C\dfrac{LB_{k+1}(\tilde{\lambda} + |p|^2)}{\sqrt{z}} \; \; .
	\end{align*}
	The estimate of the third term is then concluded.
	
	Finally, we deal with \eqref{appendix_lemma_3_second_term}, even though not necessary, we treat the cases $s>1$ and $s=1$ differently, to emphasise the influence of the exponent $2s$. Consider first $s>1$. We see that the difference in \eqref{appendix_lemma_3_second_term} is equal to
	\begin{align*}
		& \int_{\R^2}\dfrac{\hat{V}(q)(\sin \theta)^2 |q|^{2s} \ud q}{[(\tilde{\lambda} + |p|^2 + |q|^2)UB_k(\tilde{\lambda} + |p|^2 + |q|^2) + |q|^{2s}](\tilde{\lambda} + |p|^2 + |q|^2)UB_k(\tilde{\lambda} + |p|^2 + |q|^2)} \\
		& \le \int_{\R^2}\dfrac{\hat{V}(q) |q|^{2s} \ud q}{(\tilde{\lambda} + |p|^2 + |q|^2)^2UB_k(\tilde{\lambda} + |p|^2 + |q|^2)} \le C\frac{LB_k(\tilde{\lambda} + |p|^2)}{z} \int_{\R^2}\hat{V}(q) |q|^{2s-4} \ud q \le C\frac{LB_{k+1}(\tilde{\lambda} + |p|^2)}{\sqrt{z}} \; \; ,
	\end{align*}
	where in the first inequality we have used that $|q|^{2s} \ge 0$ and that $UB_k \ge 1$, and the integral in the last line is of order $\int_0^1 r^{2s-3} \ud r \le C$ since $s>1$. Now, we treat $s=1$. The difference in \eqref{appendix_lemma_3_second_term} is equal to
	\begin{align*}
		& \int_{\R^2}\dfrac{\hat{V}(q)(\sin \theta)^2 |q|^{2} \ud q}{[(\tilde{\lambda} + |p|^2 + |q|^2)UB_k(\tilde{\lambda} + |p|^2 + |q|^2) + |q|^{2}](\tilde{\lambda} + |p|^2 + |q|^2)UB_k(\tilde{\lambda} + |p|^2 + |q|^2)} \\
		& \le \int_{\R^2}\dfrac{\hat{V}(q)(\sin \theta)^2 |q|^{2} \ud q}{(\tilde{\lambda} + |p|^2 + |q|^2)^2(UB_k(\tilde{\lambda} + |p|^2 + |q|^2))^2} \le \frac{1}{\sqrt{z}} \int_{\R^2}\dfrac{\hat{V}(q)(\sin \theta)^2 \ud q}{(\tilde{\lambda} + |p|^2 + |q|^2)UB_k(\tilde{\lambda} + |p|^2 + |q|^2)} \; \; ,
	\end{align*}
	where in the last inequality we have used \eqref{appendix_CHT_exp_A1}. Now, by the estimate obtained for \eqref{appendix_lemma_3_third_term}, we have the following upper bound
	\begin{equation} \label{appendix_lemma_3_proof_second_main}
		\frac{1}{\sqrt{z}}\left[\dfrac{\pi}{2}\int_{\tilde{\lambda} + |p|^2}^1 \dfrac{\ud \varrho}{\varrho UB_k(\varrho)} + C \dfrac{LB_{k+1}(\tilde{\lambda} + |p|^2,z)}{\sqrt{z}}\right] \; \; .
	\end{equation}
	Now, note that
	\begin{equation} \label{appendix_lemma_3_squeeze_r^2_in}
		\int_{\tilde{\lambda} + |p|^2}^1 \dfrac{\ud \varrho}{\varrho UB_k(\varrho)} - \int_{\tilde{\lambda} + |p|^2}^1 \dfrac{\ud \varrho}{(\varrho^2 + \varrho)UB_k(\varrho)} = \int_{\tilde{\lambda} + |p|^2}^1 \dfrac{\ud \varrho}{(1 + \varrho) UB_k(\varrho)} \le \frac{LB_{k+1}(\tilde{\lambda} + |p|^2,z)}{\sqrt{z}}  \; \; .
	\end{equation}
	Therefore, \eqref{appendix_lemma_3_proof_second_main} is upper bounded by
	\begin{equation*} 
		\frac{1}{\sqrt{z}}\left[\dfrac{\pi}{2}\int_{\tilde{\lambda} + |p|^2}^1 \dfrac{\ud \varrho}{(\varrho^2 + \varrho)UB_k(\varrho)} + C \dfrac{LB_{k+1}(\tilde{\lambda} + |p|^2,z)}{\sqrt{z}}\right] \le C \dfrac{LB_{k+1}(\tilde{\lambda} + |p|^2,z)}{\sqrt{z}} \; \; .
	\end{equation*}
	where the last inequality is a consequence of \eqref{appendix_CHT_exp_A2} in Lemma \ref{appendix_CHT_lemma_A2}.
	The result follows from collecting all the estimates so far.
\end{proof}

\begin{lemma} \label{appendix_lemma_2}
	Let the same assumptions of Lemma \ref{appendix_CHT_lemma_A2} to hold and let $\tilde{\lambda} = \lambda + |p_{1:n}|^{2s}$, for every $s\ge1$. Then, there exists a constant $C_{\text{Off}} > 0$, such that
	\begin{equation*}
		|p| \int_{\R^2} \dfrac{\hat{V}(q)(\sin \theta)^2 \ud q}{[\tilde{\lambda} + |p+q|^2f(\tilde{\lambda} + |p+q|^2 + |q|^2)] |p^\prime + q|} \le C_{\text{Off}}\dfrac{g(\tilde{\lambda} + |p|^2)}{z} \; \; ,
	\end{equation*}
	where $p = \sum_{i=1}^{n} p_i$, $p^\prime = \sum_{i=1}^{n-1} p_i$ and $|p_{1:n}|^{2s} = \sum_{i=1}^{n} |p_i|^{2s}$.
\end{lemma}

\begin{proof}
	We split $\R^2$ into three regions, $\Omega_1 = \{q : |p + q| < \frac{|p|}{2}\}$, $\Omega_2 = \{q : |p^\prime + q| < \frac{|p|}{2}\}$ and $\Omega_3 = \R^2 \setminus (\Omega_1 \cup \Omega_2)$. Note that since we are looking for an upper bound, it is irrelevant whether $\Omega_1$ intersects $\Omega_2$ or not. In $\Omega_1$, note that $|p + q| < \frac{|p|}{2} \Rightarrow |q| < \frac{3}{2}|p| \Rightarrow |p + q|^2 + |q|^2 < \frac{5}{2}|p|^2$, so use the monotonicity of $f$ to get $f(\tilde{\lambda} + |p+q|^2 + |q|^2) \ge f(\tilde{\lambda} + \frac{5}{2}|p|^2)$ and also $(\sin \theta)^2 \le \frac{|p+q|^2}{|p|^2}$ to obtain
	\begin{align*}
		& |p| \int_{\Omega_1} \dfrac{\hat{V}(q)(\sin \theta)^2 \ud q}{[\tilde{\lambda} + |p+q|^2f(\tilde{\lambda} + |p+q|^2 + |q|^2)] |p^\prime + q|} \\
		& \le \int_{\Omega_1} \dfrac{C|p|^{-1}|p + q|^2 \ud q}{[\tilde{\lambda} + |p+q|^2f(\tilde{\lambda} + \frac{5}{2}|p|^2)] |p^\prime + q|} \le \dfrac{C |p|^{-1}}{f(\tilde{\lambda} + \frac{5}{2}|p|^2)}\int_{\Omega_1} \dfrac{\ud q}{|p^\prime + q|} \le
		C\dfrac{g(\tilde{\lambda} + |p|^2)}{z} \; ,
	\end{align*}
	since by assumption $f(x,z) \ge \frac{z}{g(x,z)}$ and $g$ is decreasing in $x$. Also, since in $\Omega_1$ we have $|p + q| < \frac{|p|}{2} \Rightarrow |q| < \frac{3|p|}{2}$, the last integral is of order $|p|$. Indeed, note that denoting $B_R(a)$ the ball of radius $R$ centred at $a$,
	\begin{equation*}
		\int_{\Omega_1} \dfrac{\ud q}{|p^\prime + q|} \le \int_{B_{\frac{3}{2}|p|}(0)} \dfrac{\ud q}{|p^\prime + q|} = 
		\int_{B_{\frac{3}{2}|p|}(-p^\prime)} \dfrac{\ud q}{| q|} \le 
		\int_{B_{6|p|}(0)} \dfrac{\ud q}{| q|} \le C \int_0^{6|p|} \ud r \le C |p| \; \; ,
	\end{equation*}
	since $q \mapsto |q|^{-1}$ has a singularity at zero.
	For the region $\Omega_2$ we use
	\begin{equation*}
		(\sin \theta)^2 \le \dfrac{4|p + q|^2}{|p|^2 \vee (\frac{1}{4}|p^\prime|^2)} \; \; ,
	\end{equation*}
	where $a \vee b := \max\{a,b\}$. This is true since, for $|p^\prime| \le 2|p|$, it is a weaker estimate than the previous one, and for $|p^\prime| > 2|p|$ it can be shown that, in the region $\Omega_2$, the right hand side is always greater or equal than $1$ (see \cite[(A.13)]{Cannizaroetal2021}). Inserting this into the integral it follows that
	\begin{equation*}
		\int_{\Omega_2} \dfrac{|p|\hat{V}(q)(\sin \theta)^2 \ud q}{[\tilde{\lambda} + |p+q|^2f(\tilde{\lambda} + |p+q|^2 + |q|^2)] |p^\prime + q|}
		\le \dfrac{C|p|}{|p|^2 \vee (\frac{1}{4}|p^\prime|^2)} \int_{\Omega_2} \dfrac{ \ud q}{f(\tilde{\lambda} + |p+q|^2 + |q|^2) |p^\prime + q|}.
	\end{equation*}
	Note that, in $\Omega_2$, we have that $|p+q|^2 + |q|^2 \le (\frac{3}{2}|p| + |p^\prime|)^2 + (\frac{1}{2}|p| + |p^\prime|)^2 \le 2(\frac{3}{2}|p| + |p^\prime|)^2$, so using the monotonicity of $f$ we obtain the upper bound
	\begin{equation} \label{appendix_proof_lemma_2_int_1}
		\dfrac{C |p|}{|p|^2 \vee (\frac{1}{4}|p^\prime|^2)f(\tilde{\lambda} + 2(\frac{3}{2}|p| + |p^\prime|)^2)} \int_{\Omega_2} \dfrac{\ud q}{|p^\prime + q|}
		= \dfrac{C |p|^2}{|p|^2 \vee (\frac{1}{4}|p^\prime|^2)f(\tilde{\lambda} + 2(\frac{3}{2}|p| + |p^\prime|)^2)} \; .
	\end{equation}
	In order to estimate the last term we maximise in $|p^\prime|$ (here we think of $p^\prime$ as any vector in $\R^2$). It can be easily seen that it is monotonously increasing for $|p^\prime| < 2|p|$. For $|p^\prime| \ge 2|p|$ we show that it is monotonously decreasing: since $f$ satisfies \eqref{appendix_CHT_exp_A7},
	for any $a,b \ge 0$, it holds that
	\begin{align*}
		& \dfrac{\ud}{\ud r} \left(\dfrac{1}{r^2f(a+ 2(b + r)^2)}\right) = -\dfrac{2rf + 4r^2(b+r)f^\prime}{r^4f^2} \\
		& = -\dfrac{2}{r^3f^2}(f + 2r(b + r)f^\prime) < -\dfrac{2}{r^3f}\left(1 - \dfrac{2r(b + r)}{a + 2(b + r)^2}\right) < 0 \; ,
	\end{align*}
	where the argument of $f$ and $f^\prime$ is always $a+ 2(b + r)^2$. Therefore, the maximum over $p^\prime$ of the right hand side of \eqref{appendix_proof_lemma_2_int_1} is attained at $|p^\prime| = 2|p|$ and is equal to
	\begin{equation*}
		\dfrac{C}{f(\tilde{\lambda} + 2(\frac{7}{2}|p|)^2)} \le C\dfrac{g(\tilde{\lambda} + |p|^2)}{z} \; \; .
	\end{equation*}
	The final part of the proof is to consider the region $\Omega_3$, for which we use $(\sin \theta)^2 \le 1$ and apply \hol inequality with exponents $\frac{3}{2}$ and $3$, to the functions $[\tilde{\lambda} + |p+q|^2f(\tilde{\lambda} + |p+q|^2 + |q|^2)]^{-1}$ and $|p^\prime + q|^{-1}$ with respect to the measure $\hat{V}(q) \ud q$ to get
	\begin{align}
		& |p| \int_{\Omega_3} \dfrac{\hat{V}(q)(\sin \theta)^2 \ud q}{[\tilde{\lambda} + |p+q|^2f(\tilde{\lambda} + |p+q|^2 + |q|^2)] |p^\prime + q|} \nonumber\\
		& \le |p| \left(\int_{\Omega_3} \dfrac{\hat{V}(q) \ud q}{[\tilde{\lambda} + |p+q|^2f(\tilde{\lambda} + |p+q|^2 + |q|^2)]^{\frac{3}{2}}}\right)^{\frac{2}{3}}
		\left(\int_{\Omega_3} \dfrac{\hat{V}(q) \ud q}{|p^\prime + q|^3}\right)^{\frac{1}{3}} \; . \label{appendix_proof_lemma_2_int_holder}
	\end{align}
	Since in $\Omega_3$ we have that $|p^\prime + q| \ge \frac{|p|}{2}$, the second term in \eqref{appendix_proof_lemma_2_int_holder} is bounded by a constant times $|p|^{-\frac{1}{3}}$.
	
	Moving to the integral inside the first parenthesis in \eqref{appendix_proof_lemma_2_int_holder}, note that in $\Omega_3$ we have that $|p + q| \ge \frac{|p|}{2} \Rightarrow |q| \le 3|p+q|$ and then by the monotonicity of $f$ we get the upper bound
	\begin{equation} \label{appendix_proof_lemma_2_int_2}
		\int_{\Omega_3} \dfrac{\hat{V}(q) \ud q}{[\tilde{\lambda} + |p+q|^2f(\tilde{\lambda} + 10|p+q|^2)]^{\frac{3}{2}}} \le \int_{\frac{|p|}{2}}^\infty \dfrac{C r \ud r}{(\tilde{\lambda} + r^2f(\tilde{\lambda} + 10r^2))^{\frac{3}{2}}}
	\end{equation}
	where the last inequality is obtained by bounding $\hat{V}(q)$ by a constant, setting $\tilde{q} = p + q$ and passing to polar coordinates. Now, we divide the domain of integration $\frac{|p|}{2} \le r < \infty$ into two regions, $\tilde{\lambda} < r^2$ and its possibly empty complement $\tilde{\lambda} \ge r^2$. In the first, it holds that
	\begin{equation} \label{appendix_proof_lemma_2_ineq_1}
		\tilde{\lambda} + r^2f(\tilde{\lambda} + 10r^2) \ge \dfrac{1}{20}(\tilde{\lambda} + 10r^2)f(\tilde{\lambda} + 10r^2) \; \; .
	\end{equation}
	Using that by assumption $f(x,z) \ge \frac{z}{g(x,z)}$ and that $g$ is decreasing in $x$, together with \eqref{appendix_proof_lemma_2_ineq_1}, we can control \eqref{appendix_proof_lemma_2_int_2} by
	\begin{align}
		& \int_{\frac{|p|}{2}}^\infty \dfrac{C r \ud r}{((\tilde{\lambda} + 10r^2)f(\tilde{\lambda} + 10r^2))^{\frac{3}{2}}}
		\le C\left(\dfrac{g(\tilde{\lambda} + \frac{5}{2}|p|^2)}{z}\right)^{\frac{3}{2}} \int_{\frac{|p|}{2}}^\infty \dfrac{r \ud r}{((\tilde{\lambda} + 10r^2))^{\frac{3}{2}}} \nonumber\\
		& \le C\left(\dfrac{g(\tilde{\lambda} + |p|^2)}{z}\right)^{\frac{3}{2}} \int_{\frac{|p|}{2}}^\infty r^{-2}\ud r \le C|p|^{-1}\left(\dfrac{g(\tilde{\lambda} + |p|^2)}{z}\right)^{\frac{3}{2}} \label{appendix_proof_lemma_2_int_3} 
	\end{align}
	The last step is to consider the region $\tilde{\lambda} \ge r^2$. For that, we have
	\begin{equation*}
		\int_{\frac{|p|}{2}}^{\sqrt{\tilde{\lambda}}} \dfrac{C r \ud r}{(\tilde{\lambda} + r^2f(\tilde{\lambda} + 10r^2))^{\frac{3}{2}}}
		\le \dfrac{1}{f(11(\tilde{\lambda} + |p|^2))^{\frac{3}{2}}} \int_{\frac{|p|}{2}}^{\infty} \dfrac{\ud r}{r^2} \le C|p|^{-1}\left(\dfrac{g(\tilde{\lambda} + |p|^2)}{z}\right)^{\frac{3}{2}} \; \; .
	\end{equation*}
	Inserting all the estimates for the region $\Omega_3$ into \eqref{appendix_proof_lemma_2_int_holder} we obtain the desired upper bound
	\begin{equation*}
		C\dfrac{g(\tilde{\lambda} + |p|^2)}{z}
	\end{equation*}
	which completes the proof.
\end{proof}

\section{End of proof of Theorem \ref{intro_theo_main_theo}} \label{section_appendix_end_proof_theo_1}

\begin{proof}[end of proof of \eqref{intro_super_dif_bounds_s>1} in Theorem \ref{intro_theo_main_theo}]
	We start with the end of the proof for the upper bound. Note that the sequence $c_{2k+1}$ in \eqref{proof_of_main_theorem_upper_bound} is monotonously decreasing and convergent, so we may replace it by its limit and merge it into the constant C below. By \eqref{appendix_CHT_exp_A9}, expression \eqref{proof_of_main_theorem_upper_bound} is bounded, by
	\begin{align}
		& Cf_{2k+1} \left( \int_\lambda^1 \dfrac{\ud \varrho}{\varrho LB_{k}(\varrho,z_{2k+1})} + \dfrac{UB_{k}(\lambda, z_{2k+1})}{\sqrt{z_{2k+1}}} \right) \nonumber \\
		& \le Cf_{2k+1} \left( \int_\lambda^1 \dfrac{\ud \varrho}{(\varrho + \varrho^2) LB_{k}(\varrho,z_{2k+1})} + \dfrac{UB_{k}(\lambda, z_{2k+1})}{\sqrt{z_{2k+1}}} \right) \nonumber \\
		& \le Cf_{2k+1} UB_k(\lambda,z_{2k+1}) \le Cf_{2k+1} \dfrac{L(\lambda, 0) + z_{2k+1}}{LB_k(\lambda,0)} \label{proof_main_theorem_upper_final}
	\end{align}
	where we have used Lemma \ref{appendix_CHT_lemma_A5} for the first inequality, \eqref{appendix_CHT_exp_A3} in Lemma \ref{appendix_CHT_lemma_A1} for the second and that $LB_k$ is increasing in $z$ for the last. Now we invoke the Central Limit Theorem applied to Poisson random variables of rate one to get
	\begin{equation*}
		\lim\limits_{k \to \infty} \sum_{i = 0}^k \dfrac{k^i}{i!}e^{-k} = \dfrac{1}{2} \; \; ,
	\end{equation*}
	which yieds that uppon the choice
	\begin{equation*}
		k = k(\lambda) = \Big\lfloor\frac{\log L(\lambda,0)}{2}\Big\rfloor
	\end{equation*}
	and recalling the definition of $LB_k$ in \eqref{anal_fock_def_LB_UB}, for $\lambda$ small enough, the bound
	\begin{equation} \label{proof_main_theorem_improvement_trick}
		\dfrac{e^{-k}}{LB_k(\lambda,0) e^{-k}} \le \dfrac{C}{\sqrt{L(\lambda,0)}} \; \; .
	\end{equation}
	Inserting the above into \eqref{proof_main_theorem_upper_final} and using the definitions of $z_{2k+1} = z_{2k+1}(1)$ and $f_{2k+1} = f_{2k+1}(1)$ in \eqref{anal_fock_def_zk_fk}, we arrive at
	\begin{equation*}
		\lambda^2 \tilde{D}(\lambda) \le C (\log L(\lambda,0))^{1+\varepsilon} \sqrt{L(\lambda,0)} \; \; ,
	\end{equation*}
	which completes the proof of the upper bound, since
	\begin{equation*}
		L(\lambda,0) = \log \left(1 + \dfrac{1}{\lambda}\right) \overset{\lambda \to 0}{\sim} |\log \lambda| \; \; .
	\end{equation*}
	Moving to the lower bound, recall \eqref{proof_main_theorem_lower_final}
	\begin{equation*}
		\dfrac{C}{f_{2k+1}} \left(LB_{k}(\lambda, z_{2k+1}) - LB_{k}(1, z_{2k+1})  - \dfrac{LB_{k}(\lambda, z_{2k+1})}{\sqrt{z_{2k+1}}} \right) 
		 \ge \dfrac{C}{f_{2k+1}} \left(LB_{k}(\lambda, z_{2k+1}) - f_{2k+1} \right)
	\end{equation*}
	where we use \eqref{anal_fock_proof_op_est_lower_ineq_LB(1,z)} and that, for $k$ large enough, $1 - \frac{1}{\sqrt{z_{2k+1}}} \ge c > 0$. Also, the $-f_{2k+1}$ term only produces a constant contribution, which can be absorved by reducing $C$ if $\lambda$ is sufficiently small. Using \eqref{proof_main_theorem_improvement_trick} with the same choice of $k$, we obtain
	\begin{equation*}
		LB_k(\lambda, 0) \ge C \sqrt{L(\lambda,0)} \; \; ,
	\end{equation*}
	which allied to the definition of $f_{2k+1} = f_{2k+1}(1)$ in \eqref{anal_fock_def_zk_fk}, concludes that
	\begin{equation*}
		\lambda^2 \tilde{D} \ge C (\log L(\lambda,0))^{-1 - \varepsilon} \sqrt{L(\lambda,0)} \; \; .
	\end{equation*}
	Therefore, \eqref{intro_super_dif_bounds_s>1} follows from \eqref{set_main_res_yaglom_rev_formula} and \eqref{set_main_res_laplace_transform_to_resolvent_eq} and the proof of Theorem \ref{intro_theo_main_theo} is concluded.
\end{proof}

\section*{Acknowledgements}
GF would like to thank B\'alint T\'oth for presenting him the paper \cite{komorowskiandolla2003} and for inspiring discussions. GF also gratefully acknowledges funding via the EPSRC Studentship 2432406 in EP/V520305/1. HW acknowledges financial support by the Royal Society through the University Research Fellowship UF140187. Moreover, GF and HW are funded by the Deutsche Forschungsgemeinschaft (DFG, German Research Foundation) under Germany’s Excellence Strategy EXC 2044 -390685587, Mathematics Münster: Dynamics–Geometry–Structure.

\bibliographystyle{abbrv}
\bibliography{documents.bib}
\end{document}